\documentclass{autart}  
\usepackage{lineno,hyperref}
\modulolinenumbers[5]

\usepackage{siunitx}
\usepackage{graphicx}
\usepackage{amssymb,amsmath}
\usepackage{lipsum}
\usepackage{mathtools}
\usepackage{cuted}
\usepackage{breqn}
\newtheorem{lemma}{Lemma}
\newtheorem{theorem}{Theorem}

\newtheorem{assumption}{Assumption}

\usepackage{lineno,hyperref}
\usepackage{graphicx}
\graphicspath{{images/}}
\begin{document}
\begin{frontmatter}

\title{Semi-globally Exponential Trajectory Tracking for a Class of Spherical Robots}

\author[mymainaddress]{T.W.U. Madhushani},
\author[mysecondaryaddress]{D.H.S. Maithripala},
\author[janakaaddress]{J.V. Wijayakulasooriya},
\author[mytertiaryaddress]{J.M. Berg}
\address[mymainaddress]{Postgraduate and Research Unit, Sri Lanka Technological Campus, CO 10500, Sri Lanka.\\  udarim@sltc.lk}
\address[mysecondaryaddress]{Dept. of Mechanical Engineering, Faculty of Engineering, University of Peradeniya, KY 20400, Sri Lanka.\\  smaithri@pdn.ac.lk}
\address[janakaaddress]{Dept. of Electrical and Electronic Engineering, Faculty of Engineering, University of Peradeniya, KY 20400, Sri Lanka.\\  jan@ee.pdn.ac.lk}
\address[mytertiaryaddress]{Dept. of Mechanical Engineering, Texas Tech University, TX 79409, USA.\\ jordan.berg@ttu.edu}

\begin{abstract}
\allowdisplaybreaks
A spherical robot consists of an externally spherical rigid body rolling on a two-dimensional surface, actuated by an auxiliary  mechanism. For a class of actuation mechanisms, we derive a controller for the geometric center of the sphere to asymptotically track any sufficiently smooth reference trajectory, with robustness to bounded, constant uncertainties in the inertial properties of the sphere and actuation mechanism, and to constant disturbance forces including, for example, from constant inclination of the rolling surface. The sphere and actuator are modeled as distinct systems, coupled by reaction forces. It is assumed that the actuator can provide three independent control torques, and that the actuator center of mass remains at a constant distance from the geometric center of the sphere. We show that a necessary and sufficient condition for such a controller to exist is that for any constant disturbance torque acting on the sphere there is a constant input such that the sphere and the actuator mechanism has a stable relative equilibrium. A geometric PID controller guarantees robust, semi-global, locally exponential stability for the position tracking error of the geometric center of the sphere, while ensuring that actuator velocities are bounded.  
\end{abstract}

\end{frontmatter}

\section{Introduction}\label{sect:Introduction}
A spherical robot consists of a main spherical component whose angular velocity can be actuated by one or more auxiliary components. For brevity we will subsequently refer to the main spherical component simply as ``the body.'' The actuating components may be internal or external to the body. Generally, motion of the body can be produced by imbalanced actuators that shift the center of gravity of the assembly, or by balanced actuators that transfer angular momentum to and from the body \cite{Crossley,Chase,Armour}. The former mechanism is referred to as \textit{barycentric actuation} and the latter mechanism is referred to as \textit{momentum actuation} \cite{Agrawal,Svinin,Borisov,Joshi,Ghariblu,Halme,Alves,Cai,Karavaev}. This paper considers the robust trajectory tracking problem for a class of spherical robots using either barycentric actuation or momentum actuation. Specifically, this paper presents a geometric PID controller to ensure that the geometric center of a spherical robot body asymptotically tracks any twice-differentiable reference trajectory, with semi-global asymptotic convergence, in the presence of constant parameter variations and disturbance forces and moments. The surface upon which the body rolls is assumed to be planar, but it is allowed to have an unknown, constant, non-zero inclination. 

A sizable body of existing results address controllability and open-loop path planning with momentum actuation \cite{Agrawal,Svinin,Borisov,Joshi,Jurdjevic,Marigo,Mukherjee}. Fewer studies consider closed-loop control or barycentric actuation. To our knowledge, the only previous controllers for barycentric actuation are the open loop path planning schemes proposed by \cite{Bloch,Cai,Karavaev} for a spherical robot rolling on a horizontal plane. To our knowledge, the only previous trajectory-tracking feedback controller is the nonlinear geometric PD momentum controller derived for an inertially symmetric spherical robot on a perfectly horizontal plane \cite{Zhao,Banavar2015}. To our knowledge, only one study \cite{Zhao} takes into account the inclination of the rolling surface. This controller combines feedback linearization with sliding mode control \cite{Zhao}. This controller of \cite{Zhao} is formulated in a single coordinate patch, and hence convergence is only guaranteed to be local. The controller of \cite{Zhao} also requires perfect knowledge of the inclination of the rolling surface. To our knowledge, the result presented in the present paper significantly extends the state of the art in feedback control of spherical robots. The class of actuation mechanisms considered here is large, although not completely general. In particular, we constrain the distance between the center of mass of the actuator and the geometric center of the body to remain constant. Despite the constraint, this actuator class includes balanced reaction wheels \cite{Svinin,Borisov,Mukherjee,Joshi,Banavar2015}, control moment gyros \cite{Halme,Ghariblu}, ``hamster-wheel'' carts driven along the inner spherical surface \cite{Karavaev}, and spherical pendulums with a fixed point at the geometric center of the body \cite{Cai}. Notably, this is the first demonstration of a tracking controller that is robust to constant inclination of the rolling surface. 

The potentially complex coupled dynamics of the actuator-body system present a challenge for controlling the spherical robot. In this paper we circumvent this obstacle by actively controlling only the location of the geometric center of the body, designated as the \textit{output system}. Here we apply geometric PID control to provide semi-global, exponential asymptotic tracking by the output system. Analogously to the linear case, geometric PID provides an intuitive and robust control framework for the control of mechanical systems \cite{MaithripalaAutomatica}. This is because the geometric representation of a mechanical system, obtained by replacing the usual time derivative by the covariant derivative, is a double integrator. Although the geometric double integrator is in general nonlinear, it preserves essential features of standard linear PID control \cite{MaithripalaAutomatica}. Geometric PID control requires a fully actuated simple mechanical system, which for brevity we refer to here as a \textit{regular} system. The complete spherical robot system of body and actuators is not regular. A previous study of tracking control for a hoop robot---a planar version of the spherical robot---introduced a procedure called \emph{feedback regularization} that allowed application of a geometric PID controller \cite{UdariACC2017}. For the 3D spherical robot considered here, the error dynamics may be \emph{split} into a linear combination of two left-invariant systems and a right-invariant system on the group of rigid body motions. These error dynamics are still not regular, due to quadratic velocity terms arising from the constraint and actuator reaction forces. Feedback regularization may be used to give this error system the form of a simple mechanical system, however in this paper we instead prove that the geometric PID controller is robustly stable to the presence of these quadratic velocity terms.

In this formulation, the system zero dynamics coincide with the controlled actuator dynamics. For practical implementation, these must remain bounded, giving the actuator characteristics a crucial role in the position tracking problem. It is shown below that for momentum actuation, boundedness of the actuator velocities cannot be guaranteed using continuous control in the presence of a non-zero constant velocity reference command, torque disturbance, or inclination of the rolling surface. On the other hand, for barycentric actuation, easily verifiable conditions ensure trajectory tracking with bounded actuation.

In Section \ref{Secn:RollingSphere} we derive the open-loop equations of motion by considering the system as comprising the rigid sphere plus each of the separate rigid bodies of the actuation mechanism. Euler's rigid body equations are used to model each of the subsystems. The exterior of the body is assumed to be perfectly spherical, but the distribution of mass on the interior need not be uniform; that is, the principal inertias need not all be equal. A \textit{no-slip} constraint is applied at the point of contact between the body and the ground; that is, the linear velocity of the contact point is assumed to be zero. Constraint forces and moments between the body, surface, and actuators are written explicitly, along with the equal and opposite reaction terms. This approach to modeling is equivalent to the Lagrange-d'Alembert principle \cite{Cai}, the constrained connection method \cite{Bullo}, and the Euler-Poincare formulation \cite{Karavaev,Bloch,BlochBook}. This chosen approach allows the spherical shell and the actuators to be treated separately, which plays a crucial role in the subsequent controller development. Additionally, this approach allows barycenter and momentum actuation to be represented in a uniform framework. 

Section \ref{Secn:IntrinsicPID} presents our main results. It begins with the tracking error dynamics for the unifying model incorporating the distinct types of actuators considered here, and concludes with the development of a geometric PID controller for trajectory tracking. Section \ref{Secn:Simus} presents simulations to demonstrate the versatility of the control framework. Subsection \ref{Secn:InnerVehicleSimu} considers a hamster-ball type actuation mechanism, Subsection \ref{Secn:GyromomentSimu} considers a balanced control moment gyro, and Subsection \ref{Secn:ReactionWheelSimu} considers a balanced reaction wheel actuator. The Appendix presents a proof of the main theorem. 

\section{Equations of Motion for the Spherical Robot} \label{Secn:RollingSphere}
We derive the equations of motion for a spherical robot with one or more barycentric or momentum actuators. Each actuator consists of a rigid body connected to the spherical body through a combination of constrained and controlled degrees of freedom. For each controlled degree of freedom, a specified force or moment may be applied by the actuator on the spherical body. For each constrained degree of freedom, the corresponding force or moment felt by the spherical body must be determined from the dynamics of the actuator system. The control and constraint forces and moments are matched by equal and opposite reaction forces and moments on the actuator body. The challenge of the spherical robot control problem largely arises from coupling caused by the actuator constraints. Here we  restrict slightly the class of actuators considered, requiring the center of mass of the actuator to remain at a constant distance from the geometric center of the sphere. 

The spherical body is assumed to be balanced -- that is, the center of mass of the body plus the body-fixed balanced momentum actuators coincides with its geometric center. The body may be inertially asymmetric -- that is, its principal moments of inertias need not be equal. The external surface of the body is assumed to be perfectly spherical, and the supporting surface upon which it rolls is assumed to be perfectly planar. The supporting surface need not be horizontal. The no-slip condition used in this paper constrains the translational velocity of the point of contact relative to the supporting surface to be zero. Our procedure can also accommodate an additional ``no-twist'' or ``rubber rolling'' constraint to require the surface normal component of the sphere angular velocity to be zero relative to the supporting surface. We do not require this constraint, however as demonstrated in the simulation section our controller is capable of enforcing this constraint.  

In this paper we represent points in space using right-handed orthonormal cartesian coordinate systems---that is, using right-handed reference frames defined by three orthogonal unit vectors. The orthonormal frame $\mathbf{e} = \{\mathbf{e}_1, \mathbf{e}_2, \mathbf{e}_3\}$ is fixed in the supporting surface, with the unit direction $\mathbf{e}_3$ coinciding with the outward surface normal. The orthonormal frame $\mathbf{b} = \{\mathbf{b}_1, \mathbf{b}_2, \mathbf{b}_3 \}$ is fixed in the spherical body, with origin coinciding with the geometric center. The orthonormal frame $\mathbf{c}_i = \{\mathbf{c}_1, \mathbf{c}_2, \mathbf{c}_3 \}$ is fixed in the $i^{\mathrm{th}}$ actuator, with origin coinciding with the actuator center of mass. We also define the $3 \times 1$ matrix $e_3\triangleq [0\, 0\, 1]^T$. 

The motion of the spherical body in space is entirely captured by specifying the motion of a specific point, $O$, of the sphere and the motion of an orthonormal frame $\mathbf{b}$ fixed to the sphere. We will choose the point $O$ to be the geometric center of the sphere and let the origin of $\mathbf{b}$ coincide with $O$. The orientation of the frame $\mathbf{b}$ is related to the earth fixed frame $\mathbf{e}$ by a rotation $R(t)$ such that $\mathbf{b}=\mathbf{e}\,R(t)$. Denote by $x$ the representation of a point $P$ in space with respect to the earth fixed frame $\mathbf{e}$ and by $X$ the representation of $P$ in the body fixed frame $\mathbf{b}$. If the representation of the position of $O$ with respect to the frame $\mathbf{e}$ at time $t$ is $o(t) \in \mathbb{R}^3$ then the Euclidean nature of space implies that $x=o+RX$. The fact that $\mathbf{e}$ and $\mathbf{b}$ are right hand oriented orthonormal frames imply that $R(t)$ is a special orthogonal matrix. Conversely, any special orthogonal matrix represents a rotation of a frame with respect to another. The Lie group $SO(3)$ represents all such possible rotations of a rigid body in space with composition of sequential rotations as the group operation. 
 
Some quantities associated with rigid body motion---including force, moment of force, angular velocity, and angular momentum---can be considered as directed line segments in space. Such quantities make sense irrespective of the choice of frame and hence are called \emph{geometric invariants}. If $\gamma$ is the representation of such an invariant in the frame $\mathbf{e}$ then we will use the notation $\Gamma$ to denote its $\mathbf{b}$ representation. These two representations are related by $\gamma = R\Gamma$. Subsequently, we refer to $\Gamma$ as the \textit{body} representation and $\gamma$ as the \textit{spatial} representation of the corresponding invariant. We shall always use lower case symbols for the latter, and upper case for the former. In this paper we denote by $R_{c_i}$ the rotation that relates the orientation of actuator $i$ with respect to the body. That is $\mathbf{c}_i = \mathbf{b} R_{c_i}(t)$ where $R_{c_i}\in SO(3)$. Since $\mathbf{b} = \mathbf{e} R(t)$ we have that $\mathbf{c}_i = \mathbf{b} R_{c_i}(t)=\mathbf{e}RR_{c_i}$ and hence define $R_i\triangleq RR_{c_i}$. 

The unit direction of gravity has the representation $-e_g$ in the frame $\mathbf{e}$. The radius, mass, inertia tensor, and angular velocity of the spherical body in body-fixed and spatial frames are, respectively, $r, m_b, \mathbb{I}_b$, $\Omega$ and $\omega$. For each actuator, the mass, inertia tensor, displacement of the actuator center of mass from the body geometric center, and angular velocity in the body-fixed and spatial frames are, respectively, $m_i$,  $\mathbb{I}_i$, $X_i$, $\Omega_i$ and $\omega_i$. Denote $l_i \triangleq \|X_i\|$, where by assumption, $l_i$ is constant. We choose the actuator-fixed frame $\mathbf{c}_i$ such that the $\mathbf{c}_3$ axis points towards the geometric center of the spherical body. Then $X_i=-l_i R_{c_i}e_3$. Because the body and actuators are rigid, $\mathbb{I}_b$ and the $\mathbb{I}_i$ are constant when written in, respectively, the body and actuator frames. These tensors are represented in the spatial frame by $\mathbb{I}_{b}^{R} \triangleq R {\mathbb{I}_b}R^T$ and $\mathbb{I}_{i}^{R_{i}} \triangleq R_i {\mathbb{I}_i}R_i^T$. The representation of the position of the center of mass of the $i^{\mathrm{th}}$ actuator is $o_i$ with respect to the the spatial frame, $o_i=o+R X_i$. The total mass of all actuators is $m_a\triangleq \sum_{i=1}^nm_i$ and $X \triangleq (\sum_{i=1}^nm_i X_i)/m_a$ is the collective center of mass of all actuators with respect to $\mathbf{b}$.

\begin{figure}[t!]
	\centering
	\begin{tabular}{c}
		\includegraphics[width=0.4\textwidth]{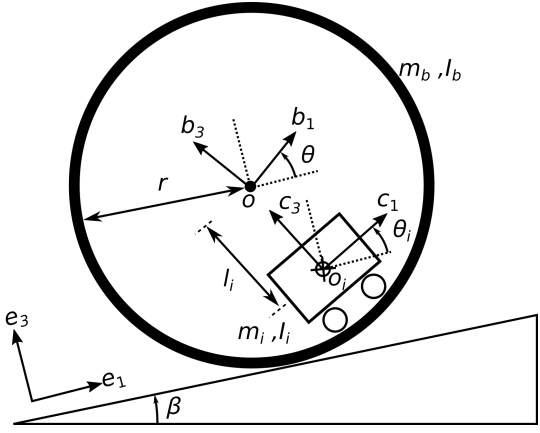}
	\end{tabular}
	\caption{Reference frames for the spherical robot analysis, for clarity shown in the plane and for a single actuator. \label{Fig:HoopRefFrames}}
\end{figure}

Actuators satisfying the constraint of constant $\|X_i\|$ include cart-driven barycentric actuators, discussed for instance in \cite{Karavaev}, barycentric pendulum actuators, discussed for instance in \cite{Cai}, and reaction-wheel momentum actuators, discussed for example in \cite{Svinin,Borisov}. Figure \ref{Fig:HoopRefFrames} depicts the problem geometry, for clarity shown in only two dimensions and for a single actuator. For barycentric actuators, the center of mass position $X_i(t)$ is non-constant. For momentum-based actuators, $X_i$ is constant. Because the system center of mass includes all fixed masses, without loss of generality we write $X_i \equiv 0$ for momentum-based actuators. 

The trajectory of the spherical body is specified by the pair $(o(t), R(t))$, and the corresponding velocity profile is $(\dot{o}(t), \dot{R}(t))$. In this paper we explicitly incorporate the constraint forces between the spherical body and the supporting surface, and so the translational velocity $\dot{o}(t)$ lies in the tangent plane to $\mathbb{R}^3$ at $o(t)$, which is also $\mathbb{R}^3$. The angular velocity $\dot{R}(t)$ lies in the tangent plane to $SO(3)$ at $R(t)$, which is isomorphic to the Lie algebra $so(3)$. The Lie algebra $so(3)$ may be identified with the vector space of $3 \times 3$ skew-symmetric matrices. It is convenient to define an isomorphism between elements of $so(3)$ and elements of $\mathbb{R}^3$. To do this, given any $X, Y \in \mathbb{R}^3$, define $\widehat{X}$ to be the unique $3 \times 3$ skew-symmetric matrix such that $\widehat{X} Y = X \times Y$ for all $Y$, where $X \times Y$ denotes the usual vector cross product. In terms of components,
\begin{equation*}
\widehat{X} =
\widehat{
\begin{bmatrix}
X_1\\
X_2\\
X_3 
\end{bmatrix}
}
=\begin{bmatrix}
0 & -X_3 & X_2\\
X_3 & 0 & -X_1\\
-X_2 & X_1 &0
\end{bmatrix}.
\end{equation*}
Then $\dot{R}_{i}=R_{i}\widehat{\Omega}_{i}=\widehat{\omega}_{i}R_{i}$, and $\dot{R}_{c_i}=R_{c_i}\widehat{\Omega}_{c_i}=\widehat{\omega}_{c_i}R_{c_i}$. Here $\Omega_i$, ${\Omega}_{c_i}$ and $\omega_i$ represent the angular velocity of the $i^{\mathrm{th}}$ actuator in, respectively, the actuator-fixed, body-fixed and spatial frames. These representations are related by ${\omega}_i={\omega}+R{\omega}_{c_i}={\omega}+R_i{\Omega}_{c_i}$.  

Given smooth vector fields $U$ and $V$, a \textit{connection} $\nabla_U V$ is a vector field that describes the change of $V$ with respect to $U$ at each point. When $U$ and $V$ are vector fields describing velocities, $\nabla_U V$ is the rate of change of $V(q)$ as the system follows the trajectories solving $\dot{q} = U$ at each point $q$. It is refered to as the \textit{covariant derivative} of $V$ along $U$. Specifically, $\nabla_U U$ can be thought of as the vector field of accelerations corresponding to velocity vector field $U$. For a Riemannian metric on a manifold, there is a unique covariant derivative referred to as the \textit{Levi-Civita connection}. For a submanifold of Euclidean space, the Levi-Civita connection on the submanifold describes the projection of the acceleration onto the tangent plane to the submanifold. Physically, this means that constraint forces corresponding to motion on the submanifold are suppressed, and do not explicitly appear in the equations of motion. 

There are many different choices of Riemannian metric for the same Lie group, which give rise to different covariant derivatives. If a metric for a Lie group is chosen so that the value of the inner product does not change under left translation, then the metric is called left-invariant. Similarly, if the metric is invariant under right translation, it is called  right-invariant. Metrics invariant under both left and right translation are called bi-invariant. In the case of rigid body motion, invariance of the metric corresponds to invariance of the inertia tensor, $\mathbb{I}_\nu$. Specifically it can be shown, using the Koszul formula \cite{Stocia}, that for left-invariant metrics induced by $\mathbb{I}_\nu$,
\[
\mathbb{I}_\nu\nabla^\nu_\xi \eta \triangleq \mathbb{I}_\nu d{\eta}(\xi)+\frac{1}{2}\left(\pm\,\mathbb{I}_\nu(\xi\times  \eta)- \left({\mathbb{I}_\nu\eta}\times {\xi} +{\mathbb{I}_\nu\xi}\times{\eta}  \right)\right),
\]
and for right-invariant metrics induced by $\mathbb{I}_\nu$,
\[
\mathbb{I}_\nu\nabla^\nu_\xi \eta \triangleq \mathbb{I}_\nu d{\eta}(\xi)+\frac{1}{2}\left(\pm\,\mathbb{I}_\nu(\xi\times  \eta)+\left({\mathbb{I}_\nu\eta}\times {\xi} +{\mathbb{I}_\nu\xi}\times{\eta}  \right)\right),
\]
and for any bi-invariant metric
\[
\nabla^{bi}_\xi \eta \triangleq d{\eta}(\xi)\pm\frac{1}{2}\,\xi\times  \eta  .
\]
When $\xi,\eta$ are body angular velocities the `$+$' are chosen; when $\xi,\eta$ are spatial angular velocities, the `$-$' are chosen. These expressions are well-defined even for singular Riemannian metrics, where $\mathbb{I}_\nu$ is positive semi-definite. For further details on the covariant derivative in general, see \cite{Marsden,Bullo,Stocia,BlochBook}. 

Subsequently we let $\nabla^b$ denote the covariant derivative on $SO(3)$ corresponding to the left-invariant Riemannian metric $\langle\langle \zeta,\eta\rangle\rangle_{so(3)}\triangleq \zeta \cdot \mathbb{I}_b\eta$ on $SO(3)$ with explicit expression  $\mathbb{I}_b^R\nabla^b_\omega \omega = \mathbb{I}_b^R\dot{\omega}-(\mathbb{I}_b^R\omega)\times \omega$. We use Newton's equations to model the body and the actuator dynamics. For the spherical body, Newton's equations specialize to Euler's equation, 
\begin{equation}
\mathbb{I}_b^R \nabla^b_\omega \omega = \tau  ,
\end{equation}
where $\tau$ is the resultant of the moments, with respect to the geometric center of the sphere, acting on the body. Euler's equation for the $i^{\mathrm{th}}$ actuator is 
\begin{equation}
\mathbb{I}_i^{R_i}\nabla^i_{{\omega_i}}{\omega_i} = \tau_i,\label{eq:Euler_ith_Actuator}
\end{equation}
where $\nabla^i$ is the unique Levi-Civita connection corresponding to the left-invariant metric $\mathbb{I}_i$ and $\tau_i$ is the resultant moment with respect to the center of mass of the actuator. We assume that each momentum-based actuator is symmetric about its actuation axis, that is, $\mathbb{I}_i$ takes the form $\mathbb{I}_i=\mathrm{diag}([\mathbb{I}_{p,i}\:\:\:\mathbb{I}_{p,i}\:\:\:\mathbb{I}_{z,i}])$. The translation of the center of mass of the rigid body obeys
\begin{equation}
m_b \ddot{o} = f,\label{eq:NewtonsCM}
\end{equation}
where $f$ is the resultant force acting on the body. The usual component-wise derivative is used here instead of the covariant derivative, because in Euclidean space the two are equivalent.


The body is acted on by the reaction forces and moments due to the actuator mechanism, forces and moments due to gravity, and the constraint forces and moments at the contact point. We denote by $f_c$ and $\tau_c$ the resultant force and moment, respectively, directly applied to the body by the actuators. We denote by $\tau_{f_c}$ the additional moment arising from the offset of the line of action of $f_c$ from the body center of mass. The gravity force acting on the sphere is $f_g=-m_b g\,e_g$ where $-e_g$ is the representation of the direction of gravity in the earth fixed frame $\mathbf{e}$.

Let the position of a point fixed on the surface of the spherical body be $y(t)$ in spatial coordinates and $Y$ in body coordinates. The velocity of such a point satisfies $\dot{y}=\dot{o}+R({\Omega}\times Y)=\dot{o}+({\omega}\times y)$, where the second equality follows because the cross product is invariant under rigid rotation. The no-slip condition requires that the velocity at the point of contact is zero, and so $\dot{o}+({\omega} \times y) = 0$. At the point of contact, $y = -re_3$, and we obtain the no-slip constraint $\dot{o} = r \widehat{e_3} \omega$. The requirement that the sphere remain in contact with the supporting surface at all times gives the constraint $e_3 \cdot \dot{o}=0$. We denote by $f_\lambda$ the constraint force acting on the body at the point of contact, and write $f=f_g + f_c + f_\lambda$. Differentiating the no-slip constraint, substituting into (\ref{eq:NewtonsCM}), and solving for $f_\lambda$ gives $f_\lambda = -\left( m_b r \widehat{e_3} \dot{\omega} + f_g + f_c \right)$. The corresponding moment about the geometric center of the body is $\tau_{f_\lambda}=-re_3\times f_\lambda = m_br^2\widehat{e_3}^2\,\dot{\omega}+re_3\times (f_g+f_c)$ and the total resultant moment acting on the body about the geometric center is $\tau=\tau_c +\tau_{f_c}+\tau_{f_\lambda}+\mathbb{I}\Delta_d$, where $\mathbb{I}\Delta_d$ represents unmodelled moments. Then the governing equation for the body is 
\begin{multline}
\mathbb{I}^R_b\nabla^b_\omega \omega-m_br^2\widehat{e_3}^2\dot{\omega} = -rm_bg\,e_3\times e_g +\mathbb{I}\Delta_d \\  + re_3 \times f_c
+\tau_c +\tau_{f_c}.   \label{eq:ConstraintBallDynamics0}
\end{multline}

The equations of motion for the actuators are similarly derived. We assume that the $i^{\mathrm{th}}$ actuator is in contact with the sphere at $p$ separate points. Let the position of the $j^{\mathrm{th}}$ contact of the $i^{\mathrm{th}}$ actuator have representation $Y_{ij}$ in the actuator frame, $\mathbf{c}_i$, with $j=1,2,\cdots,p$.  Let $-f_{ij}$ be the interaction force acting on the body at the $ij$-contact, let $-\tau_{ij}$ be the interaction moment acting at the $ij$-contact. The corresponding moment about the actuator center of mass is $-R_i Y_{ij}\times f_{ij}$, and about the geometric center of the body is $\tau_{f_{ij}}=(R X_i+R_i Y_{ij})\times f_{ij}$. Let $f_{c_i}\triangleq \sum_{j=1}^p f_{ij}$, $\tau_{c_i}\triangleq \sum_{j=1}^p \tau_{ij}$, and $\tau_{f_{c_i}}\triangleq \sum_{j=1}^p\tau_{f_{ij}}$. The total resultant moment acting about the center of mass of the actuator is $\tau_i=-\tau_{{c_i}}-\sum_{j=1}^p R_i Y_{ij}\times f_{ij}=R X_i \times f_{c_i}-(\tau_{c_i}+\tau_{f_{c_i}})$.

The total resultant reaction forces, reaction moments, and moments due to the reaction forces are, respectively, $f_{c}= \sum_{i=1}^nf_{c_i}$, $\tau_{c}= \sum_{i=1}^n\tau_{c_i}$, and $\tau_{f_{c}}= \sum_{i=1}^n\tau_{f_{c_{i}}}$. Note that all these moments are defined with respect to the geometric center of the sphere. The constraint forces $-f_{c_i}$ at the actuator contact points can be determined by twice differentiating the expression $o_i=o+R X_i$. Substituting these expressions into (\ref{eq:ConstraintBallDynamics0}) gives the final form of the constrained equations of motion for the spherical body:
\begin{multline}
\mathbb{I}_b\nabla^b_{\omega}\omega-\left((m_b+m_a)r^2\widehat{e_3}^2+rm_a\widehat{e_3} R\widehat{X}R^T\right)\dot{\omega}=\\
-rm_ae_3\times\omega\times\omega \times R X
-re_3\times\sum_{i=1}^nm_i \left(R\ddot{X}_i+2\omega\times R\dot{X}_i\right)\\
-r(m_b+m_a)g\,e_3\times e_g+\mathbb{I}\Delta_d+\sum_{i=1}^n(\tau_{c_i}+\tau_{f_{c_i}}), \label{eq:ConstraintBallDynamics}
\end{multline}
Substituting the constraint moments into (\ref{eq:Euler_ith_Actuator}) gives the equations of motion for the $i^{\mathrm{th}}$ actuator:
\begin{equation}
\mathbb{I}_i^{R_i}\nabla^i_{{\omega_i}}{\omega_i} =
R X_i \times f_{c_i} - (\tau_{c_i} + \tau_{f_{c_i}}).\label{eq:ithbodyomega}
\end{equation}

\begin{rem}
Equations (\ref{eq:ConstraintBallDynamics}) and (\ref{eq:ithbodyomega}) completely define the dynamics of an inertially asymmetric sphere rolling without slip on a, possibly inclined, planar surface. These equations are independent of the actuation mechanisms that drive the sphere. Specific types of actuators, leading to different interaction forces and moments, are considered next.
\end{rem}


{The controls act on the system through the moments $(\tau_{c_i}+\tau_{f_{c_i}})$. The exact form of the controls depend on how the actuator interacts with the sphere. For actuators interacting though rolling contact, such as cart type actuators, $f_{c_i}$ can be directly controlled and $\tau_{c_i}\equiv 0$. Hence the $\tau_{u_i}\triangleq \tau_{f_{c_i}}$ are the controls. For actuators interacting through a fixed point, such as a pendulum or a gyroscopic actuator $\tau_{c_i}$ can be directly controlled and $\tau_{f_{c_i}}\equiv 0$. For reaction wheel type actuators $\tau_{f_{c_i}}\equiv 0$ and the entire  $\tau_{c_i}$ is not available for controls. Part of the $\tau_{c_i}$ will be used to enforce the constraint that the wheel is restricted to rotate only about an axis fixed with respect to the sphere. If $\Gamma_i$ is a unit vector along this axis then we can write $\tau_{{c_i}}=u_iR\Gamma_i+\tau_{c_i}^e$ where ${u_i}\in \mathbb{R}$ is the controll and 
$-\tau_{c_i}^e$ is the moment that enforces the wheel to rotate only about $\Gamma_i$ with respect to the body.}


In the following we will consider two classes of actuators depending on whether they are balanced or not. Balanced actuators will be referred to as momentum based actuators and unbalanced actuators will be referred to as barycentric actuators. In this paper, when reaction wheels are considered, we will restrict our attention to balanced pairs of reaction wheels.  
\subsection{Barycentric Actuators}\label{Secn:BarycentricActuators} 
 
For barycentric actuators, and hence the sphere dynamics (\ref{eq:ConstraintBallDynamics}) becomes
\begin{multline}
\mathbb{I}_b^R\nabla^b_{\omega}\omega+\mathbb{I}_s\dot{\omega}
+ {e}_3\times\left(\mathbb{I}^{R_a}_{{a}}{({e}_3\times \dot{\omega})}\right) = \\
\tau_g+\tau_b+\sum_{i=1}^nB_{b_i}\tau_{u_i}, \label{eq:ConstraintBallDynamicsVehicle}
\end{multline}
where $m_a\triangleq \sum_{i=1}^n m_i$, and
$\mathbb{I}^{{R}_a}_{{a}}\triangleq\sum_{i=1}^n\mathbb{I}^{{R}_i}_{\tilde{a}_i}$,
\begin{align}
\mathbb{I}_s&\triangleq -(m_b+m_a)r^2\widehat{e_3}^2,\\
\mathbb{I}_{{a}_i}&\triangleq \mathbb{I}_{i}-m_il_i^2\widehat{e_3}^2 \\
\mathbb{I}^{{R}_i}_{\tilde{a}_i}&\triangleq-r^2m_i^2l_i^2R_i\widehat{e_3}\mathbb{I}_{a_i}^{-1}\widehat{e_3}R^T_i, \\
\tau_b &\triangleq \sum_{i=1}^n  \left(m_il_ir\,\widehat{e_3}R_i\widehat{e_3}\mathbb{I}_{a_i}^{-1}R_i^T\widehat{\omega}_i\mathbb{I}_{a_i}^{R_i}\omega_i \right. \nonumber \\
& \quad\quad + \left. m_il_ir\,\widehat{e_3}\widehat{\omega_i}^2R_i e_3 \right) \label{eq:taub}\\
\tau_g&\triangleq -r(m_b+m_a)g\,e_3\times e_g+\frac{g}{r}{e}_3\times\mathbb{I}_{{a}}^{R_a}e_g \label{eq:taug} \\
B_{b_i}&\triangleq \left(m_i l_ir\,\widehat{e_3}R_i\widehat{e_3}{\mathbb{I}_{a_i}}^{-1}R_i^T+I_{3\times 3}\right).\label{eq:Bbi}
\end{align}
For barycentric actuators the actuator dynamics (\ref{eq:ithbodyomega}) becomes 
\begin{multline}
\mathbb{I}_{{a}_i}^{R_{i}}\nabla^{a_i}_{{\omega_i}}{\omega_i}=C_{b_i}\left(\mathbb{I}_b^R\omega\times\omega+ \tau_g+\tau_b+\mathbb{I}\Delta_d\right)\\
+\tau_{g_i}- \tau_{u_i}+ C_{b_i} \sum_{k=1}^nB_{b_k}\tau_{u_k},\label{eq:ithbodyomegaVehicle0}
\end{multline}
where  $\nabla^{a_i}$ is the unique Levi-Civita connection corresponding to the left-invariant Riemannian metric induced by $\mathbb{I}_{a_i}$ and
\begin{align*}
C_{b_i}(R_i)&=-m_il_ir\,R_i\widehat{e_3}^2\left(\mathbb{I}_b^R+\mathbb{I}_{s}+\widehat{e_3}\mathbb{I}^{R_i}_{\tilde{a}_i}{\widehat{e_3}}\right)^{-1},\\
\tau_{g_i}(R_i)&=m_il_ig\,(R_i e_3)\times e_g.
\end{align*}

\subsection{Momentum Based Actuators}\label{Secn:MomentumBasedActuators}
For momentum-based actuators, since $X_i=0$, (\ref{eq:ConstraintBallDynamics}) and (\ref{eq:ithbodyomega}) become
\begin{multline}
\mathbb{I}^R_m\nabla^m_{\omega}\omega+\mathbb{I}_s\dot{\omega}=-r(m+m_a)g\,e_3\times e_g\\+\tau_{m}+\mathbb{I}\Delta_d+\sum_{i=1}^nB_{m_i}\tau_{u_i},\label{eq:ConstraintBallDynamicsCM0}
\end{multline}
\begin{equation}
A_{m_i}\mathbb{I}_{i}^{R_i}\nabla^{i}_{{\omega_i}}{\omega_i}=-\tau_{u_i},\label{eq:ithbodyomegaVehicleCM0}
\end{equation}
where $\mathbb{I}^R_m=\mathbb{I}^R_b$, $A_{m_i}=B_{m_i}=I_{3\times 3}$, and $\tau_m=0$ for non-reaction wheel type of balanced actuators and $\mathbb{I}^R_m=\mathbb{I}_r\triangleq \mathbb{I}_b-\sum_{i=1}^n\left({\widehat{\Gamma}_i^2}\mathbb{I}_i+m_i\widehat{X_i}^2\right)$, $A_{m_i}=(I_{3\times 3}+\widehat{R\Gamma}^2_i)$, $B_{m_i}=I_{3\times 3}$, and 
\begin{multline*}
\tau_m=\tau_{r}(\omega,\omega_{i}) \triangleq\\ -\sum_{i=1}^n R\widehat{\Gamma}_i^2\left({\mathbb{I}_i{R^T\omega}}\times{R^T\omega}+\mathbb{I}_i\left((R_i^T(\omega_i-\omega))\times\omega \right)\right.\\ 
\left. + \left(\mathbb{I}_i{R^T{\omega}}\right)\times{(R_i^T(\omega_i-\omega))} + (\mathbb{I}_i{R_i^T(\omega_i-\omega))}\times{R^T\omega}\right).
\end{multline*}
for balanced reaction wheel type actuators.


\section{Intrinsic Nonlinear PID Control for Trajectory Tracking}\label{Secn:IntrinsicPID}

The control problem solved in this paper is asymptotic tracking of a desired trajectory by the geometric center of the spherical body. Specifically, $o_{\mathrm{ref}}(t)$ is a twice-differentiable reference trajectory, satisfying $e_3\cdot o_{\mathrm{ref}}(t)= r$ for all $t$. This constraint has differential form $e_3\cdot \dot{o}_{\mathrm{ref}}(t)= 0$. Define error dynamics $o_e(t) \triangleq o(t) -o_{\mathrm{ref}}(t)$. Then the control objective is to ensure $\lim_{t \to \infty} o_e(t)=0$. 

For given $o_{\mathrm{ref}}(t)$ with $e_3 \cdot \dot{o}_{\mathrm{ref}} = 0$, all reference angular velocity trajectories satisfying the no-slip constraint $\dot{o}_{\mathrm{ref}} = r e_3 \times \omega_{\mathrm{ref}}(t)$ are of the form
$\omega_{\mathrm{ref}}(t;\beta) = \frac{1}{r} \left( \dot{o}_{\mathrm{ref}}(t) \times e_3 \right) + \beta e_3$
for some value of the scalar parameter $\beta$. Subsequently we assume that a smooth $\beta(t)$ has been chosen to give a suitable $\omega_{\mathrm{ref}}(t)$. For example, $\beta(t) \equiv 0$ gives the $\omega_{\mathrm{ref}}(t)$ that satisfies the ``no twist'' condition. This is the condition applied in the simulations presented below.

Now consider the control objective of ensuring that $\lim_{t \to \infty} \omega(t) = \omega_{\mathrm{ref}}(t)$. 
Let  $\omega_e\triangleq (\omega-\omega_{\mathrm{ref}})$.
Then the error dynamics can be expressed as,
\begin{align}
\dot{o}_e&=-r\,{e}_3\times\omega_e,\\
\mathbb{I}_e\dot{\omega}_e &= \mathbb{I}_b^R{{\omega}_e}\times{\omega}_e+\tau_\alpha +\tau_g+\tau_{\mathrm{ref}} +\mathbb{I}\Delta_d+ \sum_{i=1}^nB_{\alpha_i}\tau_{u_i},\label{eq:errorwe}
\end{align}
where  
\begin{align}
\mathbb{I}_e &\triangleq \left(\mathbb{I}_\alpha^R+\mathbb{I}_{s}+{e}_3\times \left(\mathbb{I}^{R_a}_{{a}}{\widehat{e_3}}\right)\right),\\
\tau_{\mathrm{ref}}&\triangleq
\mathbb{I}_\alpha^R{{\omega}_{\mathrm{ref}}}\times{{\omega}_e}+\mathbb{I}_\alpha^R{{\omega}_e}\times{{\omega}_{\mathrm{ref}}} \nonumber \\ 
&\quad\quad +\mathbb{I}_\alpha^R{{\omega}_{\mathrm{ref}}}\times{{\omega}_{\mathrm{ref}}}-\mathbb{I}_e\dot{\omega}_{\mathrm{ref}}. \label{eq:tauu}
\end{align}
\begin{rem}	
The above equations are specialized to barycentric actuators and momentum actuators by setting the subscript $\alpha=b$ and $\alpha=m$ respectively.  
\end{rem}

The objective is to find a controller that will ensure $\lim_{t\to \infty} \omega_e(t)=0$. Observe that these error dynamics do not have the structure of a mechanical system on a Lie group and thus preventing us from using the PID controller developed in \cite{MaithripalaAutomatica}.

\subsection{The zero dynamics}
The \textit{zero dynamics} of the error system are the internal dynamics of the system with $o_e$ and $\omega_e$ constrained to be identically zero. 
From (\ref{eq:ithbodyomegaVehicle0}) and (\ref{eq:ithbodyomegaVehicleCM0}) it follows that the zero dynamics of the coupled system with respect to the output $\omega_e(t)$ must satisfy
\begin{align}
A_{\alpha_i}\mathbb{I}_{{a}_i}^{R_{i}}\nabla^{a_i}_{{\omega_i}}{\omega_i}&=\tau_{g_i}({R}_i)-\tilde{\tau}_{u_i},\label{eq:ZeroDyReally}
\end{align}
where from (\ref{eq:errorwe}) we see that the output zeroing controls $\tilde{\tau}_{u_i}$ must satisfy
\begin{align}
\sum_{i=1}^nB_{\alpha_i}\tilde{\tau}_{u_i}&=-(\tau_g+\tau_\alpha+\tau_{\mathrm{ref}}+\mathbb{I}\Delta_d).\label{eq:OutputZeroControl}
\end{align}
Here we use the convention of setting the subscript $\alpha=b$ or $\alpha=m$ if the actuator is of barycentric type or momentum type respectively. Furthermore we also note that $\tau_{g_i}=0$ for momentum actuators and $A_{\alpha_i}=I_{3\times 3}$ for all actuators except reaction wheel actuators for which we have 
$A_{\alpha_i}=(I_{3\times 3}+\widehat{R\Gamma}^2_i)$.

When the output is zero, a necessary condition for the actuator states to remain bounded is that the right-hand sides of (\ref{eq:ZeroDyReally}) satisfy the condition 
$\lim_{t\to \infty}(\tau_{g_i}({R}_i)-\tilde{\tau}_{u_i})\neq \mathrm{constant}$. On the other hand it can be shown that $\lim_{t\to \infty}(\tau_{g_i}({R}_i)-\bar{\tilde{\tau}}_{u_i})= 0$ exponentially, for some constant control $\bar{\tilde{\tau}}_{u_i}$, is sufficient for the actuator states to remain bounded if there exists a positive semi-definite function $V_i : SO(3) \mapsto \mathbb{R}$ such that $dV_i=-(\tau_{g_i}(R_i)-\bar{\tilde{\tau}}_{u_i})$.

For momentum actuators, $B_{\alpha_i}=I_{3\times 3}$, $\tau_{g_i}=0$ and $\tau_g=-r(m_b+m_a)g\,e_3\times e_g$. Thus the right hand side of (\ref{eq:ZeroDyReally}) is equal to ${\tau}_{\mathrm{ref}} +\mathbb{I}\Delta_d-rg(m_b+m_a)\,\widehat{e_3}e_g$, and the necessary condition for actuator velocity boundedness is violated if the disturbances or the velocity reference are constant or if the rolling surface is not perfectly horizontal. 
\begin{rem}
For momentum actuators there exists no continuous controller that can ensure $\lim_{t\to \infty}\omega_e(t)=0$ while ensuring that the actuator velocities $\omega_i(t)$ remain bounded for non-vanishing disturbances or non-vanishing reference velocities or non-zero inclination of the rolling plane. 
\end{rem}
In particular this means that one can not stabilize the sphere at a point on an inclined plane, using a balanced actuator, while ensuring that the actuator velocities remain bounded. 
Since for momentum actuators $\tau_{g_i}=0$, the sufficient condition will be satisfied if and only if the rolling surface is perfectly horizontal and the disturbances and reference velocities tend asymptotically to zero.

We now show that for barycentric actuators, there exists a constant control $\bar{\tilde{\tau}}_{u_i}$ that ensures $\omega_e(t)\equiv 0$, and that for this $\bar{\tilde{\tau}}_{u_i}$ the actuator trajectories correspond to relative equilibria of the zero dynamics. We assume that these relative equilibria are stable---that is, there exists a positive semi-definite function $V_i : SO(3) \mapsto \mathbb{R}$ such that $dV_i=-(\tau_{g_i}(R_i)-\bar{\tilde{\tau}}_{u_i})$. Since these relative equilibria correspond to relative equilibria of rigid body rotations it is clear that they are bounded.
Thus, designing a controller that ensures $\lim_{t\to \infty}(\tau_{g_i}({R}_i)-\bar{\tilde{\tau}}_{u_i})= 0$ exponentially, guarantees  $\lim_{t\to \infty}\omega_e(t)= 0$ while ensuring that the actuator velocities $\omega_i(t)$ remain bounded.
We begin by proving the following lemma:
\begin{lemma}\label{lemma:LemmaZeroDy}
For barycentric actuators, if the disturbances, $\mathbb{I}\Delta_d$ and velocity references, $\omega_{\mathrm{ref}}$, are constant then the output satisfies ${\omega}_e(t) \equiv 0$, with ${\tau}_{u_i}(t)\equiv \mathrm{constant}$, if and only if ${\tau}_{u_i}(t)$ satisfies (\ref{eq:OutputZeroControl}) and the resulting trajectory of the actuator $({R}_i(t),{\omega}_i(t))$ satisfies ${R}_{i}(t)e_3\equiv \mathrm{constant}$. Such trajectories necessarily correspond to relative equilibria of (\ref{eq:ZeroDyReally}) with $A_{\alpha_i}=I_{3\times 3}$.
\end{lemma}

We see from (\ref{eq:taub})---(\ref{eq:Bbi}) that $B_{\alpha_i}, \tau_\alpha$ and $\tau_g$ are constant if and only if $R_{i}e_3\equiv \mathrm{constant}$. Physically what this means is that $B_{\alpha_i}, \tau_\alpha$ and $\tau_g$ are constant if and only if the actuator is  stationary with respect to the inertial frame modulo a rotation about its third body axis.
When $\omega_e\equiv 0$ equation (\ref{eq:OutputZeroControl}) implies that for constant disturbances, $\mathbb{I}\Delta_d$, constant velocity references, $\omega_{\mathrm{ref}}$, and constant control, $\tilde{\tau}_{u_i}$, the output $\omega_e\equiv 0$, if and only if $R_{i}e_3\equiv$ is constant. 
From (\ref{eq:ZeroDyReally}), with $A_{\alpha_i}=I_{3\times 3}$, we see that $(\tilde{R}_i,\tilde{\omega}_i)$ satisfies $\tilde{R}_{i}e_3\equiv \mathrm{constant}$ if and only if $(\tilde{R}_i,\tilde{\omega}_i)$ is a relative equilibrium of (\ref{eq:ZeroDyReally}).
Therefore $(\tilde{R}_i,\tilde{\omega}_i)$ must necessarily satisfy 
\begin{align}
gm_il_i(\tilde{R}_ie_3)\times e_g-\tilde{\tau}_{u_i}=0.\label{eq:EquiZeroDy}
\end{align}
Recall that (\ref{eq:OutputZeroControl}) is also necessary for the output to be zero also requires. Thus the relative equilibria $(\tilde{R}_i,\tilde{\omega}_i)$ will ensure that $\omega_e\equiv 0$ and $\tilde{R}_{i}e_3\equiv \mathrm{constant}$ if and only if $\tilde{R}_i(t)$ satisfies  
\begin{align}
&g\left(\sum_{i=1}^n\left(m_il_i\widehat{\tilde{R}_ie_3}+\frac{1}{r}\widehat{e_3}\mathbb{I}_{\tilde{a}_i}^{\tilde{R}_i}\right)-{(m_b+m_a)r}\widehat{e_3}\right)e_g \nonumber\\
&\:\:\:\:\:\:\hspace{1.5cm}\mbox{}+\tau_{\mathrm{ref}}+\mathbb{I}\Delta_d=0.\label{eq:RelativeEqm}
\end{align}

If the disturbances and velocity reference are zero, then it can be shown that there exists a $\tilde{R}_{i}(t)$ that satisfies (\ref{eq:RelativeEqm})
if the inclination $\beta\in (-\pi/2,\pi/2)$  of the inclined surface satisfies 
$\sin\beta \leq \frac{m_il_i}{(m_b+m_i)r}$,
for each actuator.  
Since for barycenter actuators $l_i < r$, this shows that for a given actuator there exists an upper bound on the inclination beyond which there exist no relative equilibria. The upper bound becomes larger for larger $l_i$ and $m_i$. In principle one can thus design an appropriate actuator such that a relative equilibrium is guaranteed to exist when the disturbances and reference velocities are constant and sufficiently small. This proves Lemma-\ref{lemma:LemmaZeroDy}.

\subsection{PID controller development}
In this section we develop a controller to ensure that the  $\lim_{t\to \infty}(o_e(t),\omega_e(t))=(0,0)$ semi-globally and exponentially while ensuring $\lim_{t\to \infty}\tilde{\tau}_{u_i}\triangleq\bar{\tilde{\tau}}_{u_i}=\mathrm{constant}$. This combined with Lemma-\ref{lemma:LemmaZeroDy} ensures the boundedness of the actuator velocities. Summarizing the conclusions of the preceding discussion, we explicitly state the conditions under which this can be ensured:
\begin{assumption}\label{as:MainAssm2}
\mbox{}
\begin{enumerate}
\item For reaction wheel actuators, the rolling plane is perfectly horizontal, and $\omega_{\mathrm{ref}}$ and $\mathbb{I}\Delta_d$ tend to zero exponentially,
\item For barycentric actuators, there exists a positive semi-definite function $V_i : SO(3) \mapsto \mathbb{R}$ such that $dV_i=-(gm_il_i({R}_ie_3)\times e_g-\bar{\tilde{\tau}}_{u_i})$ for constant $\bar{\tilde{\tau}}_{u_i}$.
\end{enumerate} 
\end{assumption}


The nonlinear PID controller, that was proposed in \cite{MaithripalaAutomatica}, for configuration tracking was based on exploiting the inherent mechanical systems structure of the error dynamics. This structure is the nonlinear equivalent of a linear double integrator. For the rolling sphere it is not possible to do this in a straightforward manner since the error dynamics given by (\ref{eq:errorwe}) do not  have the structure of an invariant mechanical system on a Lie group. 
One problem is that the candidate inertia term $\mathbb{I}_e=\left(\mathbb{I}_\alpha^R+\mathbb{I}_{s}+\widehat{e_3}\mathbb{I}^{R_a}_{{a}}{\widehat{e_3}}\right)$ is neither left-invariant nor right invariant. However we notice that it may be considered to be in some sense the sum of two left-invariant Riemannian metrics induced by $\mathbb{I}_\alpha$ and $\mathbb{I}_{{a}_i}$ and a singular right invariant metric induced by $\mathbb{I}_s$. 
Motivated by this observation if we add the quadratic velocity term
{\begin{align*}
	&\tau_{e}(\omega_e,\omega_1,\cdots,\omega_n)\triangleq \left(\mathbb{I}_{s}{\omega}_e\times{\omega}_e\right)
\nonumber\\
&\:\:\:\:\:+\frac{1}{2}\sum_{i=1}^n\left(e_3\times\left( \mathbb{I}^{R_i}_{\tilde{a}_i}({\omega_i}\times{({e}_3\times {\omega_e})})
	\right.\right.\\&\:\:\:\:\left.\left.
	+\mathbb{I}^{R_i}_{\tilde{a}_i}{\omega_i}\times{({e}_3\times {\omega_e})}
+\mathbb{I}^{R_i}_{\tilde{a}_i}{({e}_3\times {\omega_e})}\times{\omega_i}\right)\right),
	\end{align*}}
to both sides of the error dynamics (\ref{eq:errorwe}), it then takes the form
\begin{multline*}
\mathbb{I}_\alpha^R\nabla^{\alpha} _{{\omega}_e}{\omega}_e+\mathbb{I}_{s}\nabla^s_{{\omega}_e}{\omega}_e+{e}_3\times \sum_{i=1}^n\left(\mathbb{I}^{R_i}_{\tilde{a}_i}\nabla^{\tilde{a}_i}_{\omega_i}{(\widehat{e_3} {\omega_e})}\right)=\\
	\tau_e(\omega_e,\omega_1,\cdots\omega_n)+\tau_\alpha(\omega_e,\omega_1,\cdots\omega_n)\\+\tau_g +\mathbb{I}\Delta_d+\sum_{i=1}^nB_{\alpha_i}{\tau}_{u_i},
	\end{multline*}
The error dynamics written down in this form can be considered as a \textit{split mechanical system}.
Given the split mechanical structure of the error dynamics we are motivated to propose the intrinsic potential shaping plus nonlinear `split' PID controller
\begin{multline}
\mathbb{I}_s\nabla^s_{\omega_e}{o}_I+\mathbb{I}_\alpha^R\nabla^{\alpha} _{\omega_e}{o}_I\\+{e}_3\times\sum_{i=1}^n\left(\mathbb{I}_{\tilde{a}_i}^{R_i}\nabla^{\tilde{a}_i} _{\omega_i}({e}_3\times {o}_I)\right)=
\mathbb{I}_e\eta_e,\label{eq:Integrator}
\end{multline}
\begin{multline}
\sum_{i=1}^nB_{\alpha_i}{\tau}_{u_i}=\\-\left(\frac{g}{r}{e}_3\times\mathbb{I}_{{a}}^{R_a}e_{g_0}
+\mathbb{I}_e(k_p\eta_e +k_d\omega_e+k_Io_I)\right),\label{eq:PIDcontroller}
\end{multline}
where $\eta_e \triangleq \widehat{e_3}o_e$.  Here the potential shaping part of the controller given by $(\frac{g}{r}{e}_3\times\mathbb{I}_{{a}}^{R_a}e_{g_0})$ is a model of $\tau_g$ for some nominal inclination. Thus the implementation of the controller does not require the knowledge of the inclination of the rolling surface.

With this controller we see that the closed loop error dynamics are given by (\ref{eq:Integrator}) and 
{
	\begin{align}
	&\dot{o}_e=-r\,{e}_3\times \omega_e\label{eq:Erroroe}\\
	&\mathbb{I}_{s}\nabla^{s}_{{\omega}_e}{\omega}_e+\mathbb{I}_\alpha^R\nabla^{\alpha} _{{\omega}_e}{\omega}_e+\widehat{e_3}\times \sum_{i=1}^n\left(\mathbb{I}^{R_i}_{\tilde{a}_i}\nabla^{\tilde{a}_i}_{ \omega_i}{(\widehat{e_3} {\omega_e})}\right)=\nonumber\\
	&\:\:\:\:-\mathbb{I}_e(k_p\eta_e +k_d\omega_e+k_Io_I)+\mathbb{I}\Delta_d+\mathbb{I}\Delta_g +\tau_e+\tau_\alpha,\label{eq:Errorwe}
	\end{align}}
where $\mathbb{I}\Delta_g$ is given by 
\begin{align}
\mathbb{I}\Delta_g&\triangleq -r(m_b+m_a)g\,e_3\times e_g-\frac{g}{r}{e}_3\times\mathbb{I}_{{a}}^{R_a}(e_{g_0}-e_g).\label{eq:GravityDisturbance}
\end{align}
This term arises due to the ignorance of the angle of inclination in the potential shaping part of the controller $-(\frac{g}{r}{e}_3\times\mathbb{I}_{{a}}^{R_a}e_{g_0})$.
Notice that since for balanced mechanisms $\mathbb{I}_{\tilde{a}_i}=0$ this term will be absent in the controller for a balanced mechanism. However in both cases $\mathbb{I}\Delta_g\neq 0$ and hence we have the following remark:
\begin{rem}
A bounded non-vanishing unknown disturbance given by (\ref{eq:GravityDisturbance}) will always act on the error dynamics if the plane of rolling is not horizontal.
\end{rem}

In the Appendix we prove that if the PID controller gains are chosen as follows then it is possible to ensure $\lim_{t\to \infty}(o_e(t),\omega_e(t))\equiv (0,0)$ semi-globally and exponentially for bounded constant velocity references and disturbances in the presence of bounded parametric uncertainty.  
Consider a polar Morse function ${V}_\nu : \mathbb{R}^2\mapsto \mathbb{R}$ with a unique minimum at $(0,0)$ such that $dV_{\nu}=\mathbb{I}_{\nu}\eta_e$ for $\nu=\alpha,s,\tilde{a}_i$. Let $\vartheta_\nu>0$ such that
$\langle\langle\eta_e,\eta_e\rangle\rangle_{\nu}/(2\vartheta_\nu) \leq V_\nu$. Let $\vartheta_{\max}\triangleq\{\vartheta_{s},\vartheta_{b}\}$.
Let $\mu_\nu = \max \|\mathbb{I}_\nu\nabla\eta^\nu_e\|$ on $\mathcal{X}_u$,
$\mu_\min =\min \{\mu_s,\mu_b,\mu_{\tilde{a}_i}\}$, and
$\mu_\max =\max \{\mu_s,\mu_b,\mu_{\tilde{a}_i}\}$.  The existence of these constants are guaranteed since $V_\nu$ is a polar Morse function.

Let the controller gains $k_p,k_I,k_d>0$ be chosen to satisfy the following inequalities.
\begin{align}
0<k_I<\frac{k_d^3(1-\delta^2)}{\mu_{\mathrm{max}}},\label{eq:kICond}\\
 k_p>\max \left\{k_1,k_2,\frac{2 k_d^2}{\mu_{\mathrm{max}}}\right\},\label{eq:kpCond}
 \end{align}
 where, $\delta=(1-\mu_{\mathrm{min}}/\mu_{\mathrm{max}})$,
 \begin{align*}
 k_1&=\frac{k_I}{2k_d}\left(\sqrt{1+\frac{16\vartheta k_d^2}{\mu_{\mathrm{max}}^2k_I}}-1\right),\\
 k_2&=\frac{ \vartheta k_I^2}{2k_d^4}\left({1 +\sqrt{1+\frac{4k_d^3 (\mu_{\mathrm{max}}^2k_I^2+4 k_d^3(\mu_{\mathrm{max}}+k_d^3))}{\vartheta\mu_{\mathrm{max}}^2 k_I^3}}}\right).
\end{align*}

\begin{theorem}\label{theom:TheoremSphere} 
Let the conditions of Assumption-\ref{as:MainAssm2} hold and 
consider arbitrary compact sets $\mathcal{X}_s\subset \mathbb{R}^2\times  so(3)\times so(3)$ and $\mathcal{X}_a\subset SO(3)\times so(3)$. Then, if the gains of the nonlinear PID controller (\ref{eq:Integrator})--(\ref{eq:PIDcontroller}) are chosen to be sufficiently large while satisfying (\ref{eq:kICond}) and (\ref{eq:kpCond}) then for all initial conditions in $\mathcal{X}_s\times \mathcal{X}_a$, the followings hold in the presence of bounded parametric uncertainty:
\begin{enumerate}
\item if the  reference velocities and the disturbances are bounded $(o_e(t),\omega_e(t))$ converge to an arbitrarily small neighborhood of $(0,0)$,
\item if the  reference velocities and the disturbances are constant then $\lim_{t\to\infty}(o_e(t),\omega_e(t))=(0,0)$,
\end{enumerate}
while ensuring that $||\omega_i(t)||$ remains bounded. The convergence is guaranteed to be exponential.
\end{theorem}



\section{Simulation Results} \label{Secn:Simus}
In this section we simulate the performance of the intrinsic nonlinear PID controller (\ref{eq:Integrator})--(\ref{eq:PIDcontroller}). In Section \ref{Secn:InnerVehicleSimu} we simulate the performance of an inner cart actuated sphere while in Section-\ref{Secn:InnerBalancedSimu} we consider balanced actuation mechanisms.
We simulate the performance for two types of balanced mechanisms. In Section \ref{Secn:GyromomentSimu} we simulate the performance of a balanced gyroscopic moment actuated sphere and in Section \ref{Secn:ReactionWheelSimu} we simulate the performance of a balanced reaction wheel actuated sphere. 

In all simulations the nominal mass of the spherical shell was chosen to be  $m_b=1.00\, \si{kg}$, the nominal inertia tensor of the spherical shell was chosen to be $\mathbb{I}_b=\mathrm{diag}\{0.0213, 0.0205, 0.0228\} \si{kg.m^2}$, while the radius of the spherical shell was chosen to be $r=0.18\,\si{m}$. 
These parameters were chosen to correspond to a $3\, \si{mm}$ thick plastic shell, with density $850\, \si{kg.m^{-3}}$. 
To demonstrate robustness of the controller the system parameters used in the simulations were chosen to be $50\%$ different from the nominal parameters used for the controller. In all simulations the initial position of the sphere was assumed to be $o(0)=[2, -2, r]^T\,\si{m}$ and the initial angular velocity of the sphere was chosen to be $\omega(0)=[-0.1, -0.2, 0.5]^T\,\si{rad/s}$.

\subsection{Barycenter-Controlled Inner Cart Actuation}\label{Secn:InnerVehicleSimu}
Consider a  sphere actuated by an omni-directional wheel driven cart. The total mass of the cart was chosen to be $m_i=3.28\,\si{kg}$, while its inertia tensor was chosen to be
\begin{align*}
\mathbb{I}_i&=\mathrm{diag}\{0.0353, 0.0378, 0.0368\}\,\si{kg.m^2}.
\end{align*}
For these parameters one finds that the maximum inclination for which an equilibrium for the controlled cart exists is $25^\circ$.
The sphere is assumed roll on a $20^\circ$ inclined plane in the $y$-direction and is assumed to be unknown. A nominal value of $30^\circ$ inclination is assumed in the controller.

\begin{figure}[h!]
	\centering
	\begin{tabular}{cc}
		\includegraphics[width=0.235\textwidth]{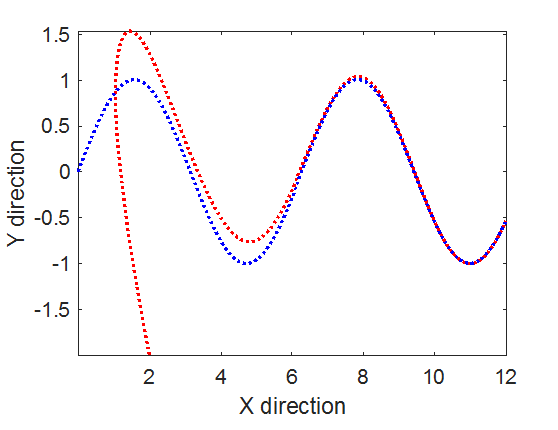}& \includegraphics[width=0.235\textwidth]{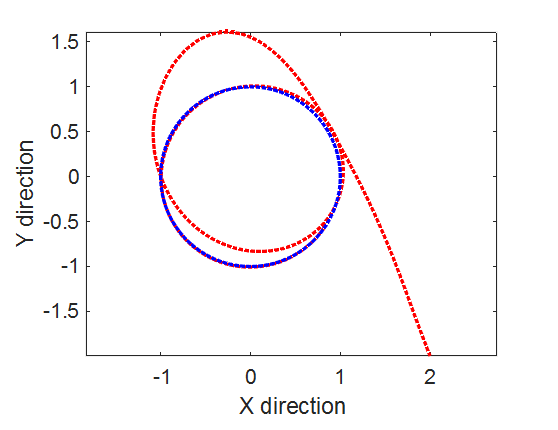}\\
		(a) Sinusoidal path & (b) Circular path\\
	\end{tabular}
	\begin{tabular}{c}
		\includegraphics[width=0.235\textwidth]{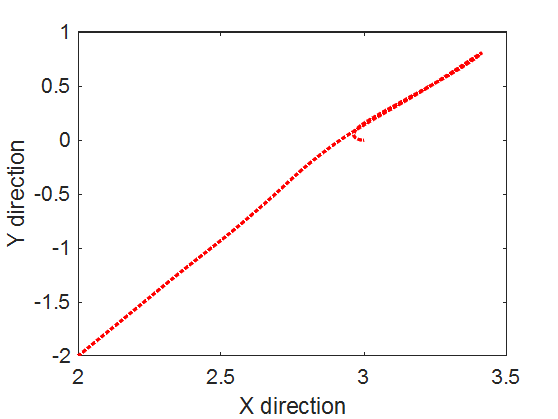}\\
		(c) Fixed point\\
		\includegraphics[width=0.235\textwidth]{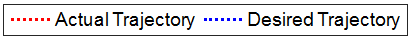}
	\end{tabular}
	\caption{The path followed by the cart actuated sphere for the PID controller (\ref{eq:Integrator})--(\ref{eq:PIDcontroller}) in the presence of parameter uncertainties as large as $50\%$. The blue curve shows the reference trajectory while the red curve shows the trajectory of the center of mass of the sphere.\label{Fig:RollingCartPath}}
\end{figure}

\begin{figure}[h!]
	\centering
	\begin{tabular}{cc}
		\includegraphics[width=0.235\textwidth]{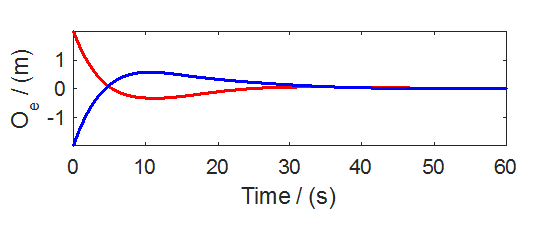}& \includegraphics[width=0.235\textwidth]{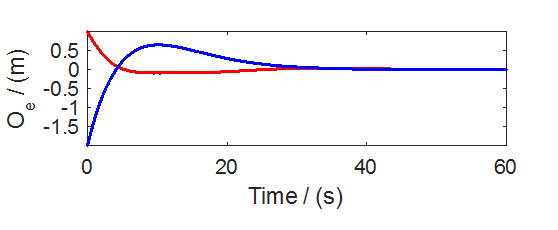}\\
		(a) Sinusoidal path & (b) Circular path\\
	\end{tabular}
	\begin{tabular}{c}
		\includegraphics[width=0.235\textwidth]{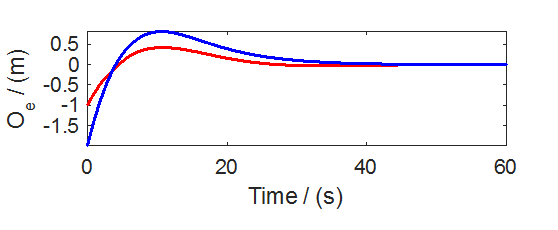}\\
		(c) Fixed point\\
		\includegraphics[width=0.1\textwidth]{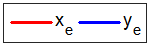}
	\end{tabular}	
	\caption{The position error, $o_e(t)$, for the cart actuated sphere for the PID controller (\ref{eq:Integrator})--(\ref{eq:PIDcontroller}) in the presence of parameter uncertainties as large as $50\%$.\label{Fig:RollingCartPositionError}}
\end{figure}

\begin{figure}[h!]
	\centering
	\begin{tabular}{cc}
		\includegraphics[width=0.235\textwidth]{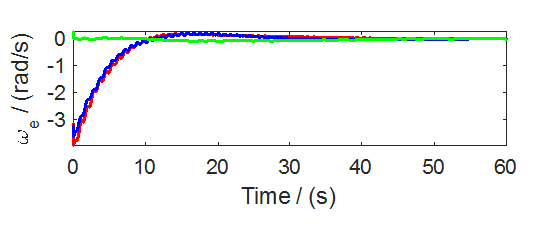}& \includegraphics[width=0.235\textwidth]{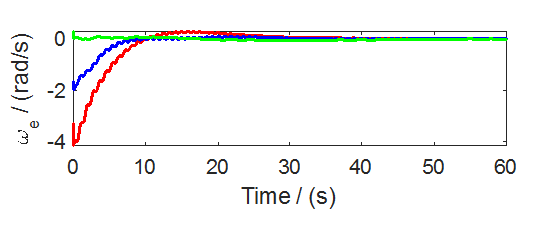}\\
		(a) Sinusoidal path & (b) Circular path\\
	\end{tabular}
	\begin{tabular}{c}
		\includegraphics[width=0.235\textwidth]{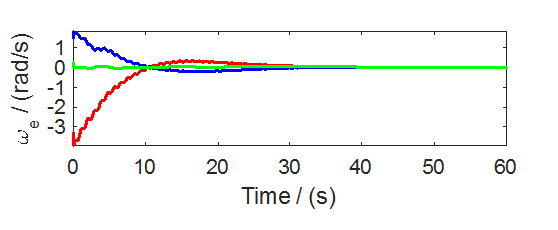}\\
		(c) Fixed point\\
		\includegraphics[width=0.1\textwidth]{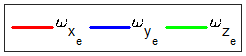}
	\end{tabular}
	\caption{The spatial angular velocity error, $\omega_e(t)$, for the cart actuated sphere for the PID controller (\ref{eq:Integrator})--(\ref{eq:PIDcontroller}) in the presence of parameter uncertainties as large as $50\%$. \label{Fig:RollingCartOmegaeError}}
\end{figure}

\begin{figure}[h!]
	\centering
	\begin{tabular}{cc}
		\includegraphics[width=0.235\textwidth]{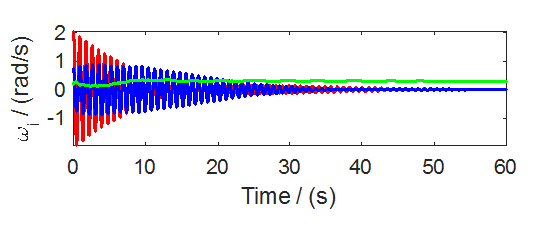}& \includegraphics[width=0.235\textwidth]{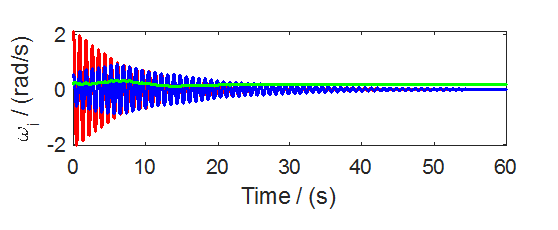}\\
		(a) Sinusoidal path & (b) Circular path\\
	\end{tabular}
	\begin{tabular}{c}
		\includegraphics[width=0.235\textwidth]{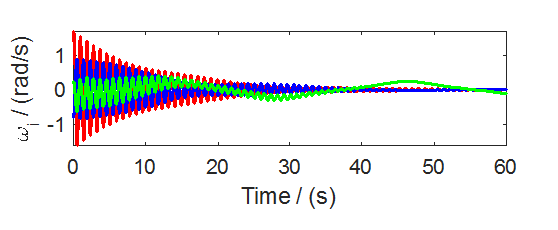}\\
		(c) Fixed point\\
		\includegraphics[width=0.1\textwidth]{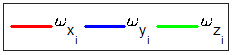}
	\end{tabular}
	\caption{The spatial angular velocities of the cart, $\omega_{i}(t)$, for the cart actuated sphere for the PID controller (\ref{eq:Integrator})--(\ref{eq:PIDcontroller}) in the presence of parameter uncertainties as large as $50\%$.\label{Fig:RollingCartomegai}}
\end{figure}

\begin{figure}[h!]
	\centering
	\begin{tabular}{cc}
		\includegraphics[width=0.235\textwidth]{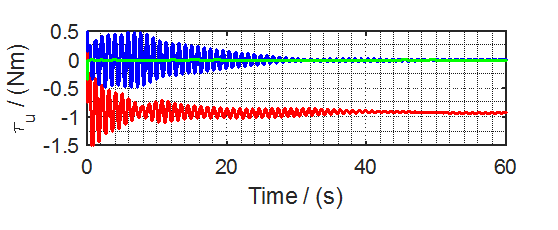}& \includegraphics[width=0.235\textwidth]{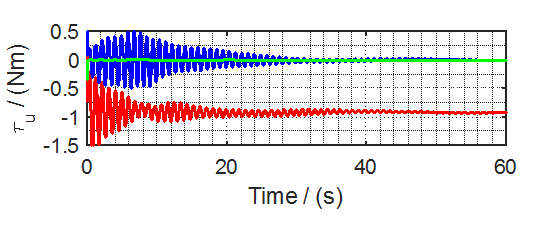}\\
		(a) Sinusoidal path & (b) Circular path\\
	\end{tabular}
	\begin{tabular}{c}
		\includegraphics[width=0.235\textwidth]{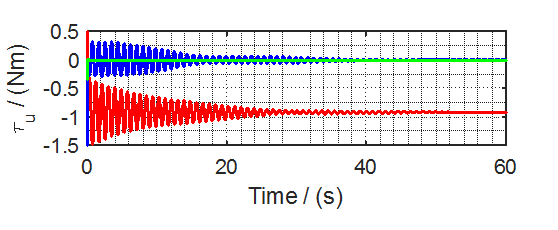}\\
		(c) Fixed point\\
		\includegraphics[width=0.1\textwidth]{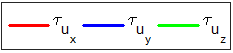}
	\end{tabular}
	\caption{Control input, $\tau_{u}(t)$, for the cart actuated sphere for the PID controller (\ref{eq:Integrator})--(\ref{eq:PIDcontroller}) in the presence of parameter uncertainties as large as $50\%$.\label{Fig:RollingCartTau}}
\end{figure}

Simulation results are presented in figure (\ref{Fig:RollingCartPath})--(\ref{Fig:RollingCartTau}) for tracking a sinusoidal path, circular path and a fixed point at $(3,0)$.
In all simulations the initial conditions used for the inner cart were $\omega_i(0)=[0.2\:\:\:-0.1\:\:\:0.1]^T\,\si{rad/s}$.
The controller gains were chosen to be $k_p=100, k_d=60, k_I=10$.  

\subsection{Balanced actuation Mechanisms}\label{Secn:InnerBalancedSimu}
In this section we simulate the behavior of the controller for two types of balanced actuation mechanisms: balanced gyroscopic moment actuation in Section-\ref{Secn:GyromomentSimu} and balanced reaction wheel actuation in Section-\ref{Secn:ReactionWheelSimu}. Simulation results are presented for tracking a sinusoidal path and a circular path. In these simulations the rolling surface is assumed to be perfectly horizontal. In both gyroscopic and reaction wheel methods the controller gains were chosen to be $k_p=55, k_d=10, k_I=1$.

\subsubsection{Balanced Gyroscopic Moment Actuation}\label{Secn:GyromomentSimu}
\begin{figure}[h!]
	\centering
	\begin{tabular}{cc}
		\includegraphics[width=0.235\textwidth]{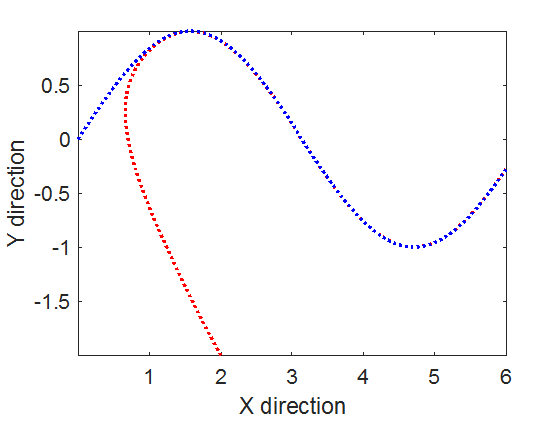}& \includegraphics[width=0.235\textwidth]{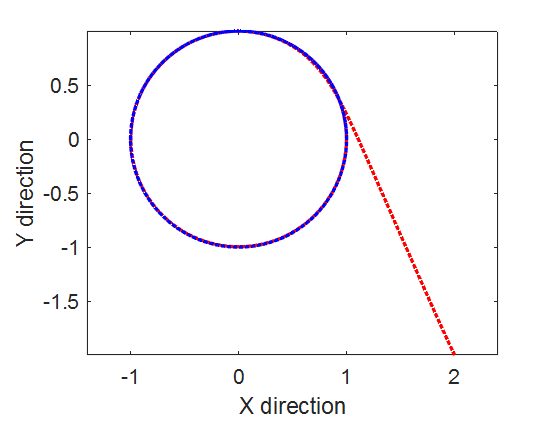}\\
		(a) Sinusoidal path & (b) Circular path
	\end{tabular}
	\begin{tabular}{c}
		\includegraphics[width=0.235\textwidth]{lgd1}
	\end{tabular}
	\caption{The path followed by the balanced gyroscopic moment actuated sphere for the PID controller (\ref{eq:Integrator})--(\ref{eq:PIDcontroller}) in the presence of parameter uncertainties as large as $50\%$. The blue curve shows the reference trajectory while the red curve shows the trajectory of the center of mass of the sphere.\label{Fig:BalancedPath}}
\end{figure}

\begin{figure}[h!]
	\centering
	\begin{tabular}{cc}
		\includegraphics[width=0.235\textwidth]{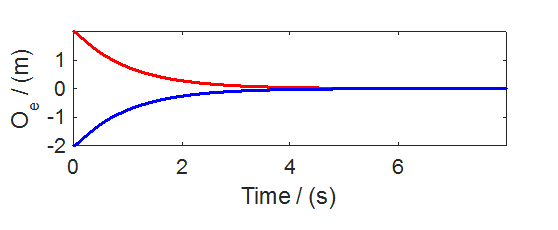}& \includegraphics[width=0.235\textwidth]{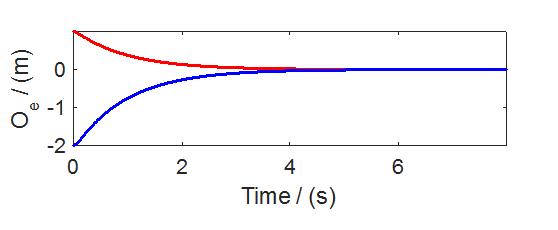}\\
		(a) Sinusoidal path & (b) Circular path
	\end{tabular}
	\begin{tabular}{c}
		\includegraphics[width=0.1\textwidth]{lgd2}
	\end{tabular}
	\caption{The position error, $o_e(t)$, for the balanced gyroscopic moment actuated sphere with the PID controller (\ref{eq:Integrator})--(\ref{eq:PIDcontroller}) in the presence of parameter uncertainties as large as $50\%$.\label{Fig:BalancedPositionError}}
\end{figure}

\begin{figure}[h!]
	\centering
	\begin{tabular}{cc}
		\includegraphics[width=0.235\textwidth]{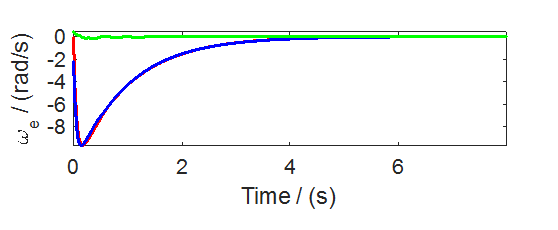}& \includegraphics[width=0.235\textwidth]{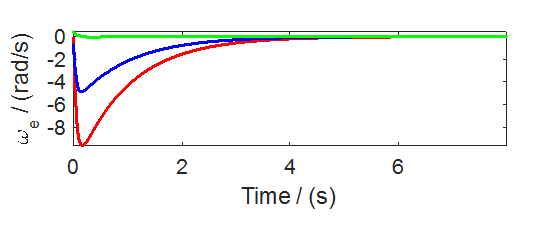}\\
		(a) Sinusoidal path & (b) Circular path
	\end{tabular}	
	\begin{tabular}{c}
		\includegraphics[width=0.1\textwidth]{lgd3}
	\end{tabular}
	\caption{The spatial angular velocity error, $\omega_e(t)$, for the balanced gyroscopic actuated sphere with the PID controller (\ref{eq:Integrator})--(\ref{eq:PIDcontroller}) in the presence of parameter uncertainties as large as $50\%$.\label{Fig:BalancedOmegaeError}
	}
\end{figure}

\begin{figure}[h!]
	\centering
	\begin{tabular}{cc}
		\includegraphics[width=0.235\textwidth]{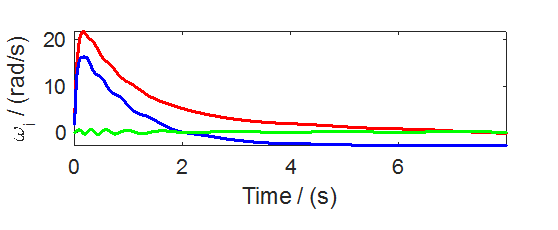}& \includegraphics[width=0.235\textwidth]{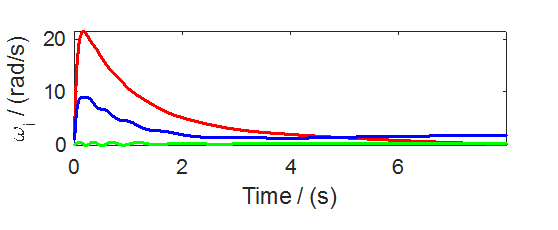}\\
		(a) Sinusoidal path & (b) Circular path
	\end{tabular}
	\begin{tabular}{c}
		\includegraphics[width=0.1\textwidth]{lgd4}
	\end{tabular}
	\caption{The spatial angular velocities of the actuation mechanism, $\omega_{i}(t)$, for the balanced gyroscopic moment actuated sphere with the PID controller (\ref{eq:Integrator})--(\ref{eq:PIDcontroller}) in the presence of parameter uncertainties as large as $50\%$. \label{Fig:GyroscopicOmegai}}
\end{figure}

In this section we present corresponding simulation results of a balanced gyroscopic moment actuator driven by omni-directional wheels similar to \cite{Halme,Ghariblu}, but with the center of mass at the geometric center of the sphere. The total mass of the mechanism was chosen to be $m_i=4.58\,\si{kg}$, while the inertia tensor was chosen to be
\begin{align*}
\mathbb{I}_i&=\mathrm{diag}\{0.0535, 0.0516, 0.0480\}\,\si{kg.m^2}
\end{align*}
In all simulations the initial conditions used for the inner vehicle were same as in section \ref{Secn:InnerVehicleSimu}. The simulation results are shown in figure (\ref{Fig:BalancedPath})--(\ref{Fig:GyroscopicOmegai}).

\subsubsection{Balanced Reaction Wheel Actuation}\label{Secn:ReactionWheelSimu}
In this section we present corresponding simulation results of actuation mechanism for balanced three pairs of reaction wheels. All pairs were not assumed to be identical. Instead we assumed the following nominal parameters for the sets of wheels: $m_1=5.78\,\si{kg}$, $m_2=4.21\,\si{kg}$, $m_3=4.91\,\si{kg}$, 
\begin{align*}
\mathbb{I}_1&=\mathrm{diag}\{0.0105, 0.0105, 0.0204\}\,\si{kg.m^2},\\
\mathbb{I}_2&=\mathrm{diag}\{0.0070, 0.0070, 0.0138\}\,\si{kg.m^2},\\
\mathbb{I}_3&=\mathrm{diag}\{0.0086, 0.0086, 0.0169\}\,\si{kg.m^2}.
\end{align*}
These parameters correspond to wheels made of lead. of radius  $r_1=0.084\,\si{m}, r_2=0.081\,\si{m}, r_3=0.083\,\si{m}$ and thickness $d_1=0.023\,\si{m}, d_2=0.018\,\si{m}, d_3=0.020\,\si{m}$ respectively.
All wheels were located at distance $l_i=0.11\, \si{m}$ from the center of the sphere. Nominal values for the rest of the components of actuation mechanism were chosen to be $m_i=4.337\,\si{kg}$,
\begin{align*}
\mathbb{I}_i&=\mathrm{diag}\{0.0166, 0.0195, 0.0053\}\,\si{kg.m^2}
\end{align*}
initial conditions used for the wheels were
$\dot{\psi}_{c_1}(0)=0.2\,\si{rad/s}$,
$\dot{\psi}_{c_2}(0)=-0.1\,\si{rad/s}$,
$\dot{\psi}_{c_3}(0)=0.1\,\si{rad/s}$.
Here we use the notation $\dot{\psi}_{c_i}=e_3\cdot\Omega_{c_i}$. 

\begin{figure}[h!]
	\centering
	\begin{tabular}{cc}
		\includegraphics[width=0.235\textwidth]{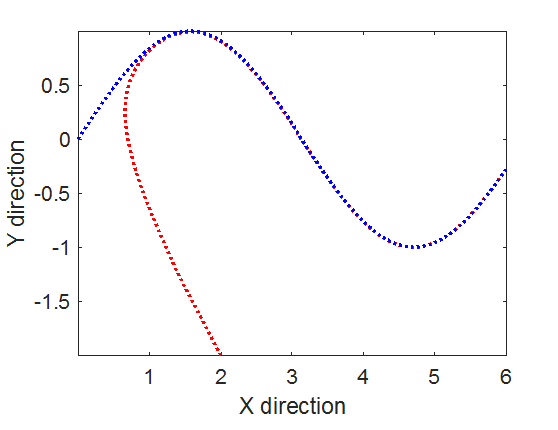}& \includegraphics[width=0.235\textwidth]{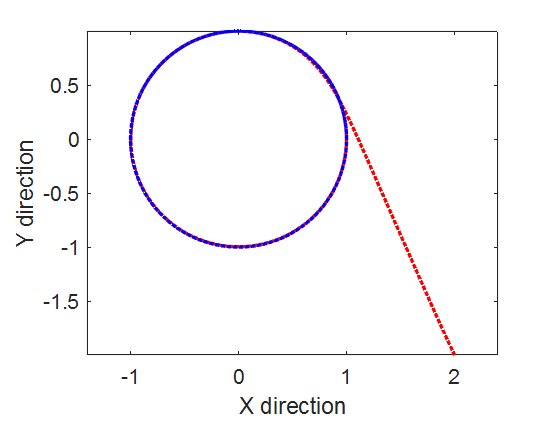}\\
		(a) Sinusoidal path & (b) Circular path
	\end{tabular}
	\begin{tabular}{c}
		\includegraphics[width=0.235\textwidth]{lgd1}
	\end{tabular}
	\caption{The path followed by the center of mass of the sphere for the reaction wheel actuated sphere with the PID controller (\ref{eq:Integrator})--(\ref{eq:PIDcontroller}) in the presence of parameter uncertainties as large as $50\%$.\label{Fig:ReactionWheelPosition}}
\end{figure}

\begin{figure}[h!]
	\centering
	\begin{tabular}{cc}
		\includegraphics[width=0.235\textwidth]{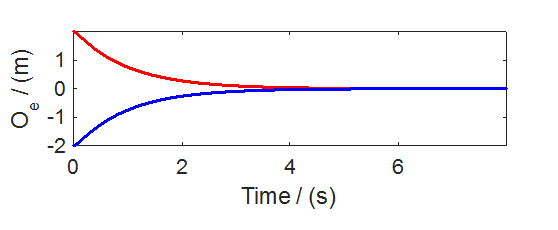}& \includegraphics[width=0.235\textwidth]{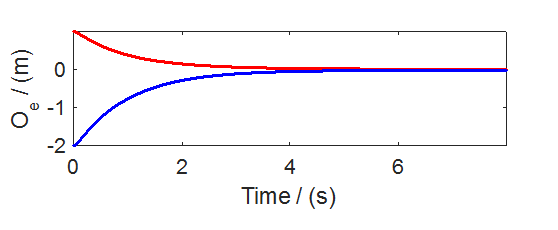}\\
		(a) Sinusoidal path & (b) Circular path
	\end{tabular}
	\begin{tabular}{c}
		\includegraphics[width=0.1\textwidth]{lgd2}
	\end{tabular}
	\caption{The position error, $o_e(t)$ of the reaction wheel actuated sphere with the PID controller (\ref{eq:Integrator})--(\ref{eq:PIDcontroller}) in the presence of parameter uncertainties as large as $50\%$.\label{Fig:ReactionWheeloe}}
\end{figure}

\begin{figure}[h!]
	\centering
	\begin{tabular}{cc}
		\includegraphics[width=0.235\textwidth]{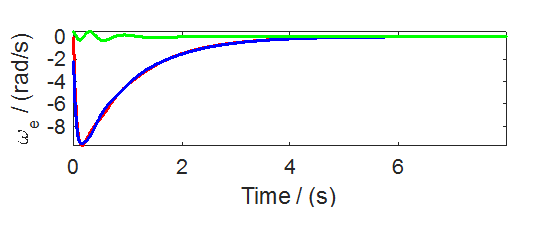}& \includegraphics[width=0.235\textwidth]{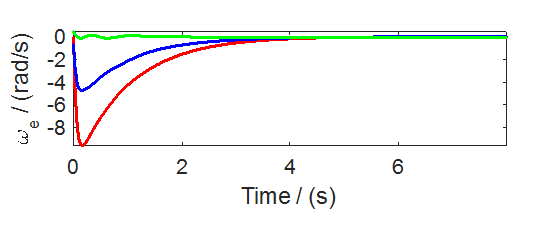}\\
		(a) Sinusoidal path & (b) Circular path
	\end{tabular}\label{Fig:ReactionWheelOmegaError}
	\begin{tabular}{c}
		\includegraphics[width=0.1\textwidth]{lgd3}
	\end{tabular}
	\caption{The angular velocity error, $\omega_e(t)$ of the reaction wheel actuated sphere with the PID controller (\ref{eq:Integrator})--(\ref{eq:PIDcontroller}) in the presence of parameter uncertainties as large as $50\%$.}
\end{figure}

\begin{figure}[h!]
	\centering
	\begin{tabular}{cc}
		\includegraphics[width=0.235\textwidth]{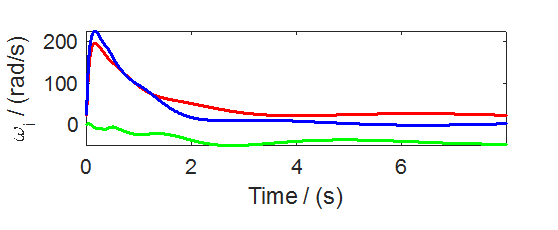}& \includegraphics[width=0.235\textwidth]{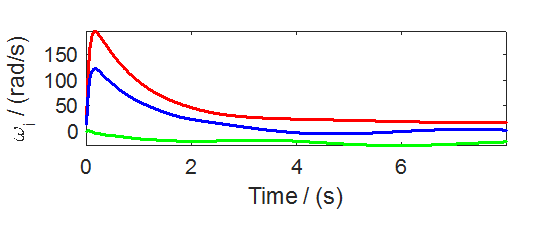}\\
		(a) Sinusoidal path & (b) Circular path
	\end{tabular}
	\begin{tabular}{c}
		\includegraphics[width=0.1\textwidth]{lgd4}
	\end{tabular}
	\caption{The angular velocities of the reaction wheels, $\Omega_{b_i}(t)$, for the reaction wheel actuated sphere with the PID controller (\ref{eq:Integrator})--(\ref{eq:PIDcontroller}) in the presence of parameter uncertainties as large as $50\%$.\label{Fig:ReactionWheelOmegai}}
\end{figure}

Figure (\ref{Fig:ReactionWheelPosition})--(\ref{Fig:ReactionWheelOmegai}) demonstrates the performance of the PID controller  (\ref{eq:Integrator})--(\ref{eq:PIDcontroller}) in the presence of parameter uncertainties as large as $50\%$. 

\section{Conclusion}
This paper considers the robust semi-global exponential tracking of an inertially non-symmetric spherical robot actuated by the class of actuation mechanisms where the distance between the center of the sphere and the center of mass of the actuation mechanism remains constant. 
The sphere is assumed to roll on a plane of unknown, constant, inclination. From a practical point of view this implies that one can robustly stabilize the sphere  at a point on an inclined surface with bounded controls using barycenter actuation. To our knowledge, this is the first appearance in the literature of a feedback controller capable of tracking a desired position on an inclined plane, in the presence of parameter uncertainty and uncertainty of the inclination of the rolling surface.

\section*{Acknowledgements}
J.\,M.\,Berg gratefully acknowledges the support of the US National Science Foundation under the Independent Research/Development program.

\bibliographystyle{elsarticle-num}
\bibliography{RollingBall2}

\begin{appendix}
\section{Proof of Theorem-\ref{theom:TheoremSphere}:}\label{append:ProofTheorem}
Consider a compact set $\mathcal{X}_s\subset \mathbb{R}^2\times so(3)\times so(3)$, constant $k_a>0$, and some small $\epsilon>0$. 
Define $\eta_e\triangleq\widehat{e_3}o_e$ and let ${V}_\nu$ be such that $dV_{\nu}=\mathbb{I}_{\nu}\eta_e$. For a given $\alpha, \beta, \gamma, \sigma >0$ define the function $W_\nu : \mathbb{R}^2\times so(3)\times so(3) \to \mathbb{R}$ as
\begin{align*}
W_\nu&\triangleq  {k_p}V_\nu(\eta_e)+\frac{1}{2}\langle\langle \omega_e,\omega_e\rangle \rangle_\nu +  \frac{\gamma}{2}\langle\langle o_I,o_I\rangle \rangle_\nu\\
&\:\:\: +\alpha \langle\langle \eta_e,\omega_e\rangle \rangle_\nu  +\beta \langle\langle o_I,\omega_e\rangle \rangle_\nu+\sigma \langle\langle o_I,\eta_e\rangle \rangle_\nu,
\end{align*}
and define the function $W : \mathbb{R}^2\times so(3)\times so(3) \to \mathbb{R}$ as
{\small
\begin{align*}
W&=W_s(\eta_e,\omega_e,o_I)+W_b(\eta_e,\omega_e,o_I)-W_{a_i}(\widehat{e_3}\eta_e,\widehat{e_3}\omega_e,\widehat{e_3}o_I).
\end{align*} 
}
where $W_s$ corresponds to the singular right invariant Reimannian metric $\langle\langle\cdot,\cdot\rangle\rangle_s$, $W_b$ corresponds to the left invariant Reimannian metric $\langle\langle\cdot,\cdot\rangle\rangle_b$, and $W_{a_i}$ corresponds to the left invariant  Reimannian metric $\langle\langle\cdot,\cdot\rangle\rangle_{a_i}$. In the following for notational simplicity we will only consider the case of a single actuator.

Let $\mathcal{W}_u$ be the set contained by the smallest level set of $W$ that contains $\mathcal{X}$
and let $k_{s}>0$ be the smallest value such that $||z_s||<k_{s}$ for all $(o_e,{\omega_e},o_I)\in \mathcal{W}_u$ where $z_s\triangleq [||\eta_s||_2\:\:\: ||\omega_e||_2\:\:\:||o_I||_2]^T$.
Let $\mathcal{W}_l$ be the set contained by the largest level set of $W$ that is contained in the set where $||z_s||<\epsilon$. 
Let $\vartheta>0$ be such that $\langle\langle\eta_e,\eta_e\rangle\rangle_{\nu}/(2\vartheta_\nu) \leq V_\nu$ on $\mathcal{W}_u$.  
Then we have 
\begin{align*}
\frac{1}{2} z_\nu^T P_l z_\nu &\leq W_\nu \leq \frac{1}{2} z_\nu^T P_u z_\nu
\end{align*}
where 
\begin{align*}
P_u=\begin{bmatrix}
\gamma &  \sigma & \beta\\
\sigma & \frac{k_p}{\vartheta_{V_i}}  & \alpha \\
\beta  &\alpha & 1
\end{bmatrix},\:\:\:\:
P_l=\begin{bmatrix}
\gamma &  -\sigma & -\beta\\
-\sigma & \frac{k_p}{\vartheta_{V_i}}  & -\alpha \\
-\beta  &-\alpha & 1
\end{bmatrix},\:\:\:\:
\end{align*}
$\vartheta_{\max}\triangleq\{\vartheta_{s},\vartheta_{b}\}$, and $z_\nu\triangleq [||o_I||_\nu\:\:\:||\eta_e||_\nu\:\:\: ||\omega_e||_\nu]^T$.
It can be shown that $P_l$ positive definite as long as 
$k_p$ satisfies (\ref{eq:kpCond}).
Since 
{\tiny
\begin{align*}
W
\geq&\frac{1}{2}\left(\left((m_b+m_i)^2r^4+\lambda^2_{\min}(\mathbb{I}_b)\right)
\lambda_{\min}(P_l)-\frac{m_i^4l_i^4r^4\lambda_{\max}(P_u)}{{\lambda^2_{\min}(\mathbb{\mathbb{I}}_i)}}
\right)||z_s||^2,
\end{align*}
}
it follows that, $W$ is positive definite,
if $k_p$ is chosen such that (\ref{eq:kpCond}) and
\[
\frac{\lambda_{\min}(\mathbb{\mathbb{I}}_i)\,\left({\lambda_{\min}(\mathbb{I}_b)+r^2(m_b+m_i)}\right)}{m_i^2l_i^2r^2}>{\frac{\lambda_{\max}(P_u)}{\lambda_{\min}(P_l)}},
\]
are satisfied.
Differentiating $W_\nu$ along the dynamics of the closed loop system we have
\begin{multline*}
\dot{W}_\nu
=\langle\mathbb{I}_\nu\nabla^\nu_{\omega_e}o_I,\beta\omega_e+\gamma o_I+\sigma \eta_e\rangle\\
 +\langle\mathbb{I}_\nu \nabla^\nu_{\omega_e}\omega_e,\omega_e+\alpha\eta_e+\beta o_I\rangle  
+ \langle \mathbb{I}_\nu\nabla^\nu_{\omega_e}\eta_e,\alpha\omega_e+\sigma o_I\rangle. 
\end{multline*}
for $\nu$ equal to $b$ or $s$ while
\begin{align*}
\dot{W}_{a_i}
&=
 \langle\mathbb{I}^{R_i}_{a_i}\nabla^{a_i}_{\omega_i}\widehat{e_3}o_I,\beta\widehat{e_3}\omega_e+\gamma \widehat{e_3}o_I+\sigma \widehat{e_3}\eta_e\rangle\\
&\:\:\: +\langle\mathbb{I}^{R_i}_{a_i} \nabla^{a_i}_{\omega_i}\widehat{e_3}\omega_e,\widehat{e_3}\omega_e+\alpha\widehat{e_3}\eta_e+\beta \widehat{e_3}o_I\rangle \\
&\:\:\: + \langle \mathbb{I}^{R_i}_{a_i}\nabla^{a_i}_{\omega_i}\widehat{e_3}\eta_e,\alpha\widehat{e_3}\omega_e+\sigma \widehat{e_3}o_I\rangle . 
\end{align*}

Then 
taking the closed loop dynamics (\ref{eq:Integrator}), (\ref{eq:Erroroe}), (\ref{eq:Errorwe}), and (\ref{eq:ZeroDyReally}) into account we can show that
\begin{align*}
\dot{W}&=\dot{W}_s+\dot{W}_\alpha-\dot{W}_{a}=DW_s+DW_b-D{W}_{a}\\
&\:\:\:\:+ \langle\epsilon_{\tau_{\mathrm{ref}}}+\epsilon\mathbb{I} (k_p\eta_e+k_d\omega_e+k_Io_I),\omega_e+\alpha\eta_e+\beta o_I\rangle,\\
&\:\:\:\:+ \langle \tau_e+\mathbb{I}\Delta_d+\mathbb{I} \Delta_g,\omega_e+\alpha\eta_e+\beta o_I\rangle,
\end{align*}
where
\begin{align*}
D{W}_\nu
&\triangleq-\beta k_I\langle \mathbb{I}_\nu o_I,o_I\rangle -(\alpha k_p-\sigma)\langle\mathbb{I}_\nu \eta_e,\eta_e\rangle\\
&\:\:\:\: -k_d\langle\mathbb{I}_\nu {\omega_e},\omega_e\rangle +(\gamma -\alpha k_I-\beta k_p)\langle\mathbb{I}_\nu o_I,\eta_e\rangle \\
&\:\:\:\:+ \langle\sigma\mathbb{I}_\nu \nabla_{\omega_e}\eta_e-(k_I+\beta k_d)\mathbb{I}_\nu\omega_e,o_I,\rangle \\
&\:\:\:\:+ \alpha\langle\mathbb{I}_\nu\nabla_{\omega_e}\eta_e,\omega_e\rangle +(\beta -\alpha k_d)\langle\mathbb{I}_\nu \omega_e,\eta_e\rangle,
\end{align*}
for $\nu$ equal to $b$ and $s$ while
\begin{align*}
&D{W}_{a_i}
\triangleq-\beta k_I\langle \mathbb{I}^{R_i}_{a_i}\widehat{e_3} o_I,\widehat{e_3}o_I\rangle -(\alpha k_p-\sigma)\langle\mathbb{I}^{R_i}_{a_i} \widehat{e_3}\eta_e,\widehat{e_3}\eta_e\rangle\\
&\:\:\:\: -k_d\langle\mathbb{I}_{a_i} \widehat{e_3}{\omega_e},\widehat{e_3}\omega_e\rangle +(\gamma -\alpha k_I-\beta k_p)\langle\mathbb{I}^{R_i}_{a_i} \widehat{e_3}o_I,\widehat{e_3}\eta_e\rangle \\
&\:\:\:\:+ \langle\sigma\mathbb{I}^{R_i}_{a_i} \nabla^{a_i}_{\omega_i}\widehat{e_3}\eta_e-(k_I+\beta k_d)\mathbb{I}^{R_i}_{a}\widehat{e_3}\omega_e,\widehat{e_3}o_I,\rangle \\
&\:\:\:\:+ \alpha\langle\mathbb{I}^{R_i}_{a_i}\nabla_{\omega_i}\widehat{e_3}\eta_e,\widehat{e_3}\omega_e\rangle +(\beta -\alpha k_d)\langle\mathbb{I}^{R_i}_{a_i} \widehat{e_3}\omega_e,\widehat{e_3}\eta_e\rangle.
\end{align*}
Let $\mu_\nu = \max ||\mathbb{I}_\nu\nabla\eta^\nu_e||$ on $\mathcal{X}_u$,
$\mu_\min =\min \{\mu_s,\mu_b,\mu_{a_i}\}$, and
$\mu_\max =\max \{\mu_s,\mu_b,\mu_{a_i}\}$.
Also let
\begin{align*}
\alpha &=\frac{k_I}{k_d^2},\:\:\:\:
\beta = \frac{k_I}{k_d},\:\:\:\:
\gamma= \frac{k_I(k_I+k_pk_d)}{k_d^2}, \:\:\:\: \sigma=\frac{2k_I}{\mu}.
\end{align*}
Then using Lemma-\ref{lemma:le1} proven below, and the properties of the inner product we have
$-z^T_\nu Q_uz_\nu\leq  D{W}_\nu \leq -z^T_\nu Q_lz_\nu$,
for all $(o_e,{\omega_e},o_I)\in \mathcal{W}_u$
where for $\delta \triangleq (1-\mu_{\min}/\mu_{\max})$ 
 \begin{align*}
 Q_l&=\begin{bmatrix}
 \frac{k_I^2}{k_d} & 0 & -\delta k_I\\
 0 &   \left(\alpha k_p-2k_I/\mu_{\min} \right) & \frac{(k_I-\alpha k_d^2)}{2k_d} \\
-\delta k_I & \frac{(k_I-\alpha k_d^2)}{2k_d} & k_d- {\alpha\mu_{\max}}
 \end{bmatrix},\\
 Q_u&=\begin{bmatrix}
 \frac{k_I^2}{k_d} & 0 & \delta k_I\\
 0 &   \left(\alpha k_p-2k_I/\mu_{\max} \right) & -\frac{(k_I-\alpha k_d^2)}{2k_d} \\
\delta k_I & -\frac{(k_I-\alpha k_d^2)}{2k_d} & k_d- {\alpha\mu_{\min}}
 \end{bmatrix}.
 \end{align*}
It can be shown that  we can pick gains such that $\lambda_{\min}(Q_\nu)$ is arbitrary and $\alpha,\beta <1$.

Since $\tau_e$ is quadratic in the velocity we see that
\[
\langle \tau_e,{\omega_e}+\alpha{\eta_s}+\beta {o_I}\rangle \leq g_1||\omega_i||^2||z_s||+g_2||\omega_i||||z_s||^2+g_3||z_s||^3.
\]
Then we have
{\tiny
\begin{align*}
\dot{W}
& \leq -\left(\lambda_{\mathrm{min}}(Q_l)-9g_0\epsilon_{\mathbb{I}}\right)(||z_s||_b^2+||z_s||_s^2)+\left(\lambda_{\mathrm{min}}(Q_u)+9g_0\epsilon_{\mathbb{I}}\right)||z_s||_{a}^2\\
&\:\:\:\:+\epsilon_{\mathrm{ref}}||z_s||^3+g_1||\omega_i||^2_2||z_s||+g_2||\omega_i||_2||z_s||^2+g_3||z_s||^3\\
&\:\:\:\:+3(||\mathbb{I}\Delta_d||+||\mathbb{I}\Delta_g||)||z_s||\\
& \leq -\left(\chi_l(\lambda^2_{\min}(\mathbb{I}_b)+r^4(m_b+m_i)^2)-\frac{\chi_ur^4l^4_im_i^4}{{\lambda^2_{\min}(\mathbb{I}_i)}}-g_2\right)||z_s||^2\\
&\:\:\:\:+(\epsilon_{\mathrm{ref}}+g_3)||z_s||^3+\left(3(||\mathbb{I}\Delta_d||+||\mathbb{I}\Delta_g||)+g_1||\omega_i||^2_2\right)||z_s||
\end{align*}}
for all $(o_e,{\omega_e},o_I)\in \mathcal{W}_u$ where 
$\chi_l\triangleq\left(\lambda_{\mathrm{min}}(Q_l)-g_0\epsilon_{\mathbb{I}}\right)$, 
$\chi_u\triangleq\left(\lambda_{\mathrm{max}}(Q_u)+g_0\epsilon_{\mathbb{I}}\right)$, and $\epsilon_{\mathbb{I}}$ and $\epsilon_{\mathrm{ref}}$ are small constants that depend on the uncertainty of the knowledge of the system parameters.

We will show below that there exists a compact set $\mathcal{X}_a\subset SO(3)\times so(3)$ such that $\lim_{t\to \infty}(o_e(t),\omega_e(t),o_I(t))=(0,0,\bar{o}_I)$ exponentially and $||\omega_i(t)||_2<k_a$ for all $t>0$ for any  $(o_e(0),\omega_e(0),o_I(0))\in \mathcal{X}_s$ and $(R_i(0),\omega_i(0))\in \mathcal{X}_a$. 
Let $g_4\triangleq \left(3(||\mathbb{I}\Delta_d||+||\mathbb{I}\Delta_g||)+g_1k_{a}^2\right)$, and
{\tiny\begin{align*}
\Xi&\triangleq\left(\chi_l\lambda^2_{\min}(\mathbb{I}_b)+r^4\left(\chi_l(m_b+m_i)^2-\frac{\chi_um_i^4l_i^4r^4}{\lambda^2_{\min}(\mathbb{I}_i)}\right)-g_2-(g_3+\epsilon_{\mathrm{ref}})k_{s}\right).
\end{align*}}
Then we have $\dot{W} \leq  -(\Xi||z_s||-g_4)||z_s||$ on $\mathcal{W}_u$ provided that $||\omega_i(t)||<k_a$.
The right hand side of this inequality is negative if
{\[
||z_s||>\epsilon_c=\frac{\left(3(||\mathbb{I}\Delta_d||+||\mathbb{I}\Delta_g||)+g_1k_{a}^2\right)}{\Xi}.
\]}
If $\lambda_{\mathrm{min}}(Q_l),\lambda_{\mathrm{min}}(Q_u)>0$ the
condition 
\begin{align*}
\frac{\lambda_{\min}(\mathbb{\mathbb{I}}_i)\,\left({\lambda_{\min}(\mathbb{I}_b)+r^2(m_b+m_i)}\right)}{m_i^2l_i^2r^2}>{\frac{\lambda_{\mathrm{min}}(Q_u)+g_0\epsilon_{\mathbb{I}}}{\lambda_{\mathrm{min}}(Q_l)-g_0\epsilon_{\mathbb{I}}}}
\end{align*} 
ensures that $\Xi>0$. We have shown that it is possible to find gains $k_p,k_I,k_d$ such that $\lambda_{\mathrm{min}}(Q_l),\lambda_{\mathrm{min}}(Q_u)$ are arbitrarily large. 
Thus for a given $\epsilon ,k_{s},k_a>0$ it is possible to find gains $k_p,k_I,k_d$ such that $\epsilon_c<\epsilon $. Thus it is possible to find gains $k_p,k_I,k_d$ such that
$\dot{W}<0$ on $\mathcal{W}_u/\mathcal{W}_l$ for bounded disturbances and reference velocities and parametric uncertainty, provided that $||\omega_i(t)||<k_a$ and the actuator satisfies
\begin{align*}
&\frac{\lambda_{\min}(\mathbb{\mathbb{I}}_i)\,\left({\lambda_{\min}(\mathbb{I}_b)+r^2(m_b+m_i)}\right)}{m_i^2l_i^2r^2}\\
&>\max\left\{{\frac{\lambda_{\max}(P_u)}{\lambda_{\min}(P_l)}},{\frac{\lambda_{\mathrm{min}}(Q_u)+g_0\epsilon_{\mathbb{I}}}{\lambda_{\mathrm{min}}(Q_l)-g_0\epsilon_{\mathbb{I}}}}\right\}
\end{align*}
Without loss of generality we will assume that at the design stage the actuator is choses such that ${\lambda_{\min}(\mathbb{\mathbb{I}}_i)}$ is sufficiently larger than ${m_il_i^2}$.
Thus proving that $(o_e(t),\omega_e(t),o_I(t))$ can be made to converge to an arbitrarily small neighborhood of $(0,0,0)$ provided that $||\omega_i(t)||<k_a$. Since $W$ is quadratically bounded from below the convergence is guaranteed to be exponential.

If the disturbances and the velocity references are constant then the Lasalle's invariance theorem says that the trajectories converge to the largest invariant set in the set where $\dot{W}\equiv 0$.
These invariant sets are exactly the equilibrium solutions of (\ref{eq:Erroroe}) -- (\ref{eq:Errorwe}) and (\ref{eq:Integrator})  and take the form $(0,0,\bar{o}_I)$. Thus we see that 
if the disturbances and the velocity references are constant then $\lim_{t \to \infty}(o_e(t),\omega_e(t),o_I(t))=(0,0,\bar{o}_I)$ for all $(o_e(0),\omega_e(0),o_I(0))\in\mathcal{X}$ provided that $||\omega_i(t)||<k_a$. The convergence is exponential since $W$ is bounded below by a positive definite quadratic form.

We will now show that the exponential convergence of 
$\lim_{t\to \infty} (o_e(t),{\omega_e}(t),{o_I}(t))=(0,0,\bar{o}_I)$ guarantees the existence of a $k_a>0$ such that $||\omega_i(t)||\leq k_a$ for all $t>0$ and $(R_i(0),\omega_i(0))\in \mathcal{X}_a$.
The exponential convergence implies that there exists $\kappa>0$ such that $||\omega_e(t)||\leq ||\omega_e(0)||e^{-\kappa t}$. Since from Lemma-\ref{lemma:LemmaZeroDy} and (\ref{eq:PIDcontroller}) it follows that $\lim_{t\to \infty}\omega_e(t)=0$ implies that $\lim_{t\to \infty}(\tilde{\tau}_{u_i}(t)-\bar{\tilde{\tau}})\triangleq \tau_a=0$ exponentially, we have that there exists $\nu>0$ such that $||\tau_a||<\nu e^{-\kappa t}$. This constant $
\kappa$ is an increasing function of $\left(\lambda_{\mathrm{min}}(Q_l)-g_0\epsilon_{\mathbb{I}}\right)$ and since we have shown that $\lambda_{\mathrm{min}}(Q_l)$ can be arbitrarily assigned it follows that we can make $\kappa$ arbitrarily large as well provided that the parametric uncertainty is sufficiently small. 


Consider the positive semi-definite function $W_i: SO(3)\times so(3)\mapsto \mathbb{R}$ defined to be $W_i=V_i+{\langle \langle \omega_i,\omega_i\rangle\rangle}/2$.
The derivative of this function along the dynamics of the closed loop system satisfies
$\dot{W}_i\leq 2\left(\nu +k_sg_i\right)e^{-\kappa t}\,W_i$.
This gives that 
$W_i\leq W_i(0)\,\exp{\left(\frac{2\left(\nu+k_sg_i\right)}{\kappa}\right)}$.
Hence we have that 
\[
||v_a(t)||\leq {(2V_a(0)+||v_a(0)||^2)}\,\exp{\left(\frac{\left(\nu+k_sg^a_1\right)}{\kappa}\right)}.
\] 
Let $k_a^0>0$ be such that ${(2V_i(R_i)+||\omega_i||^2)}\leq k_a^0$ on $\mathcal{X}_a$.
We have shown that by picking sufficiently large PID gains $k_p,k_I,k_d$ one can make $\lambda_{\min}(Q_s)$ and hence $\kappa$ sufficiently large. Thus there exists gains such that 
${(2V_i(0)+||\omega_i(0)||^2)}\,e^{\frac{\left(\nu+k_sg_i\right)}{\kappa}}$ can be made less than $k_a>k_a^0$ for all $(R_i,\omega_i)\in \mathcal{X}_i$ and hence ensure that $||\omega_i(t)||<k_a$ for all time $t>0$. 
\qed

\begin{lemma}\label{lemma:le1}
Let $V$ be a Morse function on a smooth Riemannian manifold $G$ with a metric induced by the inertia tensor $\mathbb{I}$. Let $\eta\in TG$ be the gradient of $V$ (ie. $dV=\mathbb{I}\eta$). For a given compact $\mathcal{K}\subset {G}$ there exists $\kappa >0$ and $0\leq \delta \leq1$ such that $||\kappa\mathbb{I}\nabla\eta-I_{n\times n}||\leq \delta$ on $\mathcal{K}$.
\end{lemma}

\textit{Proof of Lemma-\ref{lemma:le1}}\\
\\
Since
\begin{align*}
\langle\langle \zeta,\kappa\nabla_{\zeta}\eta-\zeta\rangle \rangle  & = \kappa\langle\langle \zeta,\nabla_{\zeta}\eta \rangle \rangle-||\zeta||^2 \leq(\kappa ||\mathbb{I}\nabla \eta||-1)||\zeta||^2,
\end{align*}
we have that 
$||\kappa\nabla\eta-I_{n\times n}||=|(\kappa ||\mathbb{I}\nabla \eta||-1)| < 1$
as long as $0 <\kappa ||\mathbb{I}\nabla \eta || < 2$.
By assumption $V$ is a Morse function and hence it is non-degenerate at all its critical points thus we have that $ \mu_{\mathrm{min}}<||\mathbb{I}\nabla\eta|| < \mu_{\mathrm{max}}$ for some $0\leq\mu_{\mathrm{min}}<\mu_{\mathrm{max}}$ on $\mathcal{K}$.
This implies that $\kappa$ should satisfy
$0 <\kappa  < {2}/{\mu_{\mathrm{max}}}$.
Choose $\kappa  = {1}/{\mu_{\mathrm{max}}}$. Then we have that
\begin{align*}
||\kappa\nabla\eta-I_{n\times n}||&=|(\kappa ||\mathbb{I}\nabla \eta||-1)| \leq \delta\triangleq \left(1-\frac{\mu_{\mathrm{min}}}{\mu_{\mathrm{max}}}\right).
\end{align*}
on $\mathcal{K}$. \qed


\end{appendix}
\end{document}